\documentclass[a4,12pt]{article}
\usepackage{amsmath,amsfonts,amssymb,amsthm,latexsym,diagrams}
\usepackage{amscd,graphicx,stmaryrd}
\begin{document}
\baselineskip15.25pt
\newcommand{\nc}[2]{\newcommand{#1}{#2}}
\newcommand{\rnc}[2]{\renewcommand{#1}{#2}}
\rnc{\theequation}{\thesection.\arabic{equation}}
%\rnc{\section}{\setcounter{equation}{0}}
\def\note#1{{}}
%\def\Label#1{\label{#1}\ifmmode\llap{[#1] }\else
%\marginpar{\smash{\hbox{[#1]}}}\fi}

%%%%%%%%%%%%%%%%%%% THEOREMS %%%%%%%%%%%%%%%%%% 
\newtheorem{definition}{Definition $\!\!$}[section]
\newtheorem{proposition}[definition]{Proposition $\!\!$}
\newtheorem{lemma}[definition]{Lemma $\!\!$}
\newtheorem{corollary}[definition]{Corollary $\!\!$}
\newtheorem{theorem}[definition]{Theorem $\!\!$}
\newtheorem{example}[definition]{\sc Example $\!\!$}
\newtheorem{remark}[definition]{\sc Remark $\!\!$}
%%%%%%%%%%%%%%%%%%%%%%%%%%%%%%%%%%%%%%%%%%%% 

%%%%%%%%%%%%%% ENVIROMENT %%%%%%%%%%%%%%%%%%%%%%
\nc{\beq}{\begin{equation}}
\nc{\eeq}{\end{equation}}
\rnc{\[}{\beq}
\rnc{\]}{\eeq}
\nc{\ba}{\begin{array}}
\nc{\ea}{\end{array}}
\nc{\bea}{\begin{eqnarray}}
\nc{\beas}{\begin{eqnarray*}}
\nc{\eeas}{\end{eqnarray*}}
\nc{\eea}{\end{eqnarray}}
\nc{\be}{\begin{enumerate}}
\nc{\ee}{\end{enumerate}}
\nc{\bd}{\begin{diagram}}
\nc{\ed}{\end{diagram}}
\nc{\bi}{\begin{itemize}}
\nc{\ei}{\end{itemize}}
\nc{\bpr}{\begin{proposition}}
\nc{\bth}{\begin{theorem}}
\nc{\ble}{\begin{lemma}}
\nc{\bco}{\begin{corollary}}
\nc{\bre}{\begin{remark}\em}
\nc{\bex}{\begin{example}\em}
\nc{\bde}{\begin{definition}}
\nc{\ede}{\end{definition}}
\nc{\epr}{\end{proposition}}
\nc{\ethe}{\end{theorem}}
\nc{\ele}{\end{lemma}}
\nc{\eco}{\end{corollary}}
\nc{\ere}{\hfill\mbox{$\Diamond$}\end{remark} }
\nc{\eex}{\hfill\mbox{$\Diamond$}\end{example}}
\nc{\bpf}{{\it Proof.~~}}
\nc{\epf}{\hfill\mbox{$\square$}\vspace*{3mm}}
\nc{\hsp}{\hspace*}
\nc{\vsp}{\vspace*}
\newcommand{\wegdamit}[1]{}
%%%%%%%%%%%%%%%%%%%%%%%%%%%%%%%%%%%%%%

%%%%%%%%%%%%%%%%% SYMBOLS %%%%%%%%%%%%%%%
%\def\span{{\rm span}}
\nc{\ot}{\otimes}
\nc{\te}{\!\ot\!}
\nc{\bmlp}{\mbox{\boldmath$\left(\right.$}}
\nc{\bmrp}{\mbox{\boldmath$\left.\right)$}}
\nc{\LAblp}{\mbox{\LARGE\boldmath$($}}
\nc{\LAbrp}{\mbox{\LARGE\boldmath$)$}}
\nc{\Lblp}{\mbox{\Large\boldmath$($}}
\nc{\Lbrp}{\mbox{\Large\boldmath$)$}}
\nc{\lblp}{\mbox{\large\boldmath$($}}
\nc{\lbrp}{\mbox{\large\boldmath$)$}}
\nc{\blp}{\mbox{\boldmath$($}}
\nc{\brp}{\mbox{\boldmath$)$}}
\nc{\LAlp}{\mbox{\LARGE $($}}
\nc{\LArp}{\mbox{\LARGE $)$}}
\nc{\Llp}{\mbox{\Large $($}}
\nc{\Lrp}{\mbox{\Large $)$}}
\nc{\llp}{\mbox{\large $($}}
\nc{\lrp}{\mbox{\large $)$}}
\nc{\lbc}{\mbox{\Large\boldmath$,$}}
\nc{\lc}{\mbox{\Large$,$}}
\nc{\Lall}{\mbox{\Large$\forall\;$}}
\nc{\bc}{\mbox{\boldmath$,$}}
\nc{\ra}{\rightarrow}
\nc{\ci}{\circ}
\nc{\cc}{\!\ci\!}
\nc{\lra}{\longrightarrow}
%\rnc{\to}{\mapsto}
\nc{\imp}{\Rightarrow}
\rnc{\iff}{\Leftrightarrow}
\nc{\inc}{\mbox{$\,\subseteq\;$}}
\rnc{\subset}{\inc}
\def\sw#1{{\sb{(#1)}}}
\nc{\0}{\sb{(0)}}
\nc{\1}{\sb{(1)}}
\nc{\2}{\sb{(2)}}
\nc{\3}{\sb{(3)}}
\nc{\4}{\sb{(4)}}
\nc{\5}{\sb{(5)}}
\nc{\6}{\sb{(6)}}
\nc{\7}{\sb{(7)}}
\newcommand{\Boxneu}{\square}
\def\tr{{\rm tr}}
\def\Tr{{\rm Tr}}
\def\st{\stackrel}
\def\<{\langle}
\def\>{\rangle}
\def\d{\mbox{$\mathop{\mbox{\rm d}}$}}
\def\id{\mbox{$\mathop{\mbox{\rm id}}$}}
\def\ker{\mbox{$\mathop{\mbox{\rm Ker$\,$}}$}}
\def\coker{\mbox{$\mathop{\mbox{\rm Coker$\,$}}$}}
\def\hom{\mbox{$\mathop{\mbox{\rm Hom}}$}}
\def\im{\mbox{$\mathop{\mbox{\rm Im}}$}}
\def\map{\mbox{$\mathop{\mbox{\rm Map}}$}}
\def\o{\sp{[1]}}
\def\t{\sp{[2]}}
\def\mo{\sp{[-1]}}
\def\z{\sp{[0]}}
\def\lhom#1#2#3{{{}\sb{#1}{\rm Hom}(#2,#3)}}
\def\rhom#1#2#3{{{\rm Hom}\sb{#1}(#2,#3)}}
\def\lend#1#2{{{}\sb{#1}{\rm End}(#2)}}
\def\rend#1#2{{{\rm End}\sb{#1}(#2)}}
\def\Rhom#1#2#3{{{\rm Hom}\sp{#1}(#2,#3)}}
\def\Lhom#1#2#3{{{}\sp{#1}{\rm Hom}(#2,#3)}}
\def\Rrhom#1#2#3#4{{{\rm Hom}\sp{#1}\sb{#2}(#3,#4)}}
\def\Llhom#1#2#3#4{{{}\sp{#1}\sb{#2}{\rm Hom}(#3,#4)}}
\def\LRhom#1#2#3#4#5#6{{{}\sp{#1}\sb{#2}{\rm Hom}\sp{#3}\sb{#4}(#5,#6)}}
\def\khom#1#2{{{\rm Hom}(#1,#2)}}
%%%%%%%%%%%%%%%%%%%%%%%%%%%%%%%%%%%%%%%%%%%%

%%%%%%%%%%%%%%%%%% OBJECTS %%%%%%%%%%%%%%%%%%%%
\nc{\spp}{\mbox{${\cal S}{\cal P}(P)$}}
\nc{\ob}{\mbox{$\Omega\sp{1}\! (\! B)$}}
\nc{\op}{\mbox{$\Omega\sp{1}\! (\! P)$}}
\nc{\oa}{\mbox{$\Omega\sp{1}\! (\! A)$}}
\nc{\dr}{\mbox{$\Delta_{R}$}}
\nc{\dsr}{\mbox{$\Delta_{\Omega^1P}$}}
\nc{\ad}{\mbox{$\mathop{\mbox{\rm Ad}}_R$}}
\nc{\as}{\mbox{$A(S^3\sb s)$}}
\nc{\bs}{\mbox{$A(S^2\sb s)$}}
\nc{\slc}{\mbox{$A(SL(2,\C))$}}
\nc{\suq}{\mbox{$\cO(SU_q(2))$}}
\nc{\tc}{\widetilde{can}}
\def\slq{\mbox{$\cO(SL_q(2))$}}
\def\asq{\mbox{$\cO(S_{q,s}^2)$}}
\def\esl{{\mbox{$E\sb{\frak s\frak l (2,{\Bbb C})}$}}}
\def\esu{{\mbox{$E\sb{\frak s\frak u(2)}$}}}
\def\ox{{\mbox{$\Omega\sp 1\sb{\frak M}X$}}}
\def\oxh{{\mbox{$\Omega\sp 1\sb{\frak M-hor}X$}}}
\def\oxs{{\mbox{$\Omega\sp 1\sb{\frak M-shor}X$}}}
\def\Fr{\mbox{Fr}}
\nc{\p}{{\rm pr}}
\newcommand{\can}{\operatorname{\it can}}
%%%%%%%%%%%%%%%%%%%%%%%%%%%%%%%%%%%%%%%%%%%%

%%%%%%%%%%%%%%%% HELLENIC %%%%%%%%%%%%%%%%%%%%%%
%\rnc{\varepsilonilon}{\varepsilon}
\rnc{\phi}{\varphi}
\nc{\ha}{\mbox{$\alpha$}}
\nc{\hb}{\mbox{$\beta$}}
\nc{\hg}{\mbox{$\gamma$}}
\nc{\hd}{\mbox{$\delta$}}
\nc{\he}{\mbox{$\varepsilon$}}
\nc{\hz}{\mbox{$\zeta$}}
\nc{\hs}{\mbox{$\sigma$}}
\nc{\hk}{\mbox{$\kappa$}}
\nc{\hm}{\mbox{$\mu$}}
\nc{\hn}{\mbox{$\nu$}}
\nc{\hl}{\mbox{$\lambda$}}
\nc{\hG}{\mbox{$\Gamma$}}
\nc{\hD}{\mbox{$\Delta$}}
\nc{\hT}{\mbox{$\Theta$}}
\nc{\ho}{\mbox{$\omega$}}
\nc{\hO}{\mbox{$\Omega$}}
\nc{\hp}{\mbox{$\pi$}}
\nc{\hP}{\mbox{$\Pi$}}
%%%%%%%%%%%%%%%%%%%%%%%%%%%%%%%%%%%%%%%%%%%%

%%%%%%%%%%%%%%%%% PHRASES %%%%%%%%%%%%%%%%%%%%%
\nc{\qpb}{quantum principal bundle}
\def\gal{-Galois extension}
\def\hge{Hopf-Galois extension}
\def\ses{short exact sequence}
\def\csa{$C^*$-algebra}
\def\ncg{noncommutative geometry}
\def\wrt{with respect to}
\def\Ha{Hopf algebra}
%%%%%%%%%%%%%%%%%%%%%%%%%%%%%%%%%%%%%%%%%%%%

%%%%%%%%%%%%%%%% FONTS %%%%%%%%%%%%%%%%%%%%%%%%
\def\C{{\Bbb C}}
\def\N{{\Bbb N}}
\def\R{{\Bbb R}}
\def\Z{{\Bbb Z}}
\def\T{{\Bbb T}}
\def\Q{{\Bbb Q}}
\def\cO{{\mathcal O}}
\def\cT{{\cal T}}
\def\cA{{\cal A}}
\def\cD{{\cal D}}
\def\cB{{\cal B}}
\def\cK{{\cal K}}
\def\cH{{\cal H}}
\def\cM{{\cal M}}
\def\cJ{{\cal J}}
\def\fB{{\frak B}}
\def\pr{{\rm pr}}
\def\ta{\tilde a}
\def\tb{\tilde b}
\def\td{\tilde d}
%%%%%%%%%%%%%%%%%%%%%%%%%%%%%%%%%%%%%%%%%%%%

\title{\LARGE\bf \vspace*{-15mm}
THE K-THEORY OF HEEGAARD-TYPE QUANTUM 3-SPHERES\\ 
%\footnote{PRELIMINARY VERSION}
\vspace*{3mm}{\Large\it Dedicated to the memory of Olaf Richter.}}
\author{
\vspace*{-0mm}\Large\sc
Paul Baum   \\
\vspace*{-0mm}\large
Mathematics Department,
McAllister Building  \\
\vspace*{-0mm}\large
The Pennsylvania State University  \\
\vspace*{-0mm}\large
University Park, PA  16802, USA   \\
\vspace*{-0mm}\large\sl e-mail: baum@math.psu.edu  \\
%\large http://www.math.psu.edu/baum\\
\vspace*{0mm}\and
\vspace*{-0mm}\Large\sc
Piotr M.~Hajac\\
\vspace*{-0mm}\large
Instytut Matematyczny, Polska Akademia Nauk\\
\vspace*{-0mm}\large
ul.\ \'Sniadeckich 8, Warszawa, 00-956 Poland\\
\vspace*{-0mm}\large
and\\
\vspace*{-0mm}\large
Katedra Metod Matematycznych Fizyki, Uniwersytet Warszawski\\
\vspace*{-0mm}\large
ul.\ Ho\.za 74, Warszawa, 00-682 Poland \vspace{0mm}\\
\large\sl
http://www.fuw.edu.pl/$\!\widetilde{\phantom{m}}\!$pmh\\
\vspace*{0mm}\and
\vspace*{-0mm}\Large\sc
Rainer Matthes\\
\vspace*{-0mm}\large
Katedra Metod Matematycznych Fizyki, Uniwersytet Warszawski\\
\vspace*{-0mm}\large
ul.\ Ho\.za 74, Warszawa, 00-682 Poland \vspace{0mm}\\
\large\sl
e-mail: matthes@itp.uni-leipzig.de\\
\vspace*{0mm}\and
\Large\sc
Wojciech Szyma\'nski\\
\vspace*{-0mm}\large
School of Mathematical and Physical Sciences, University of Newcastle\\
\vspace*{-0mm}\large
Callaghan, NSW 2308, Australia\\
\large\sl
e-mail: wojciech@frey.newcastle.edu.au
}
\date{\normalsize }

\maketitle

{\large
\begin{abstract}\normalsize
We use a Heegaard splitting of the topological 3-sphere as a guiding principle
to construct a family of its noncommutative deformations. The main technical point is
an identification of the universal $C^*$-algebras defining our quantum 3-spheres
with an appropriate fiber product of crossed-product $C^*$-algebras.
Then we employ this result to show that the $K$-groups  of
our family of noncommutative 3-spheres coincide with their classical counterparts. 
\end{abstract}
}

\section*{Introduction}
\setcounter{equation}{0}

It seems we are in a period of development of Noncommutative Geometry when 
 new constructions of quantum spheres appear in abundance. Perhaps a key reason
 for this richness of examples is that a wish to preserve some particular additional
 structures that live on a topological space often leads to, or allows for, significantly
 different types of noncommutative deformations. The flexibility and scope of the language
 of noncommutative $C^*$-algebras leaves enough room to reflect much more in the 
 structure of a $C^*$-algebra than just pure topology of the underlying space.

 For instance,  the concept of
 the $SU(2)$
  group structure on the 3-sphere uniquely determined (under some natural assumptions)
   its deformation into the quantum group $SU_q(2)$ 
 \cite{w-sl87}. Among recent
 examples, we have the classification of quantum deformations of $S^3$ that preserve
 certain classical cohomological conditions \cite{cd02,cd03}. 
 Both new and earlier examples of noncommutative spheres are beautifully surveyed
in \cite{d-l03}, and we refer the reader for more details and references therein.
 
 Our present work stems from the idea of the decomposition of $S^3$ into two solid tori that
 corresponds to the local triviality of the Hopf fibration $S^3\ra S^2$. 
 It is known as a Heegaard splitting of $S^3$. This decomposition
 principle was first explored in \cite{m-k91} yielding the quantum-torus type 
 deformation $S^3_\theta$, and, recently, in \cite{cm02,hms} producing the quantum-disc
 type deformation $S^3_{pq}$. Both approaches were unified  into
 $S^3_{pq\theta}$ in \cite{bhms}, where the coordinate algebra $\cO(S^3_{pq\theta})$
  is studied from Hopf-Galois and index  theory
point of view. Herein, we compute the $K$-theory of the universal enveloping $C^*$-algebra
  of $\cO(S^3_{pq\theta})$, which is the main result of this paper. 

To understand better the $C^*$-algebraic structure we work with, note first
 that the aforementioned Heegaard splitting of $S^3$ decomposes it into the glueing
over the boundary torus of two copies of a
solid torus which is the Cartesian product of a disc and 
 circle. Our noncommutative
3-sphere is an analogous glueing (represented by the fiber product of algebras) of
two copies of a quantum solid torus. We represent a quantum solid torus by the
crossed product of the Toeplitz algebra $\cT$ by a natural action of $\Z$ determined
by the parameter $\theta$, i.e.,  $\cT {\rtimes}_ {\theta}\,\Z$. Just as the standard
epimorphism $\cT\ra C(S^1)$ allows one to think of $S^1$ as the boundary of a
quantum disc, the natural epimorphism 
$\cT {\rtimes}_ {\theta}\,\Z\ra C(S^1){\rtimes}_ {\theta}\,\Z$
permits one to view the noncommutative torus as the boundary of a quantum solid
torus. This is the way the Toeplitz and irrational rotation $C^*$-algebras are used
to assemble a new $C^*$-algebra along the lines of
 the topological idea of a Heegaard splitting of~$S^3$.

Before embarking on a computation of the $K$-groups of noncommutative deformations
of topological spaces, it is reasonable to ask under which general conditions we have
a guarantee that the $K$-groups of the deformed $C^*$-algebras remain the same as
the $K$-groups of the initial $C^*$-algebras. An example of such general 
conditions was provided
in \cite{r-ma93}. On the other hand, even a very standard deformation of a torus into
a quantum torus $T^2_\theta$, after passing to the quotient by the antipodal action
$T^2_\theta/Z_2$, leads to different $K$-theory \cite{beek91,beek92,bk92}: 
\[
K_0(C(T^2_\theta/Z_2))=K_0((C(S^1)\rtimes_\theta\Z)^{Z_2})\cong\Z^6
\neq\Z^2\cong K^0(S^2)= K_0(C(T^2/Z_2)).
\]
Therefore, it is in general not true that the $K$-theory remains unchanged with deformations.
As it seems unlikely that our example fits into Rieffel's framework \cite{r-ma93}, we prove
directly that the $K$-groups of our family of quantum spheres coincide with their classical
counterparts.

The overall mathematical structure of this article is much in line with \cite{m-k91}.
Section-wise, it is organized as follows. We begin with preliminaries where notation is fixed
and basic facts recalled. Also in preliminaries, we warm up for the $K$-theory computation
by analysing the classical geometric roots of the Mayer-Vietoris six-term
 exact sequence in the $K$-theory
of $C^*$-algebras. Then we proceed with the main technical contents of the paper, which
is the determination of key $C^*$-algebraic features of our quantum spheres. First we show
that the universal $C^*$-algebra remains unchanged if we set two of our parameters to zero,
i.e., $C(S^3_{pq\theta})\cong C(S^3_{00\theta})$. Next we prove the latter $C^*$-algebra
to be isomorphic with a fiber product (pullback) of two crossed products. These two steps
make it possible to compute the $K$-groups of $C(S^3_{pq\theta})$ by standard
methods of the $K$-theory of operator algebras \cite{w-ne93,b-b98}.

\section{Preliminaries}
\setcounter{equation}{0}

Throughout the paper we use the jargon of noncommutative geometry referring
to quantum spaces as objects dual  to noncommutative algebras
in the sense of the Gelfand-Naimark correspondence
between spaces and function algebras. 
 The algebras are assumed to be associative and over $\C$.
The notation
$C_0(\mbox{locally compact Hausdorff space})$ means the 
algebra of vanishing-at-infinity continuous functions on this space. 
This algebra is non-unital unless the space is compact, in which case we skip
the zero subscript. (In the noncommutative setting, the unitality of a $C^*$-algebra
is viewed as the compactness of a quantum space.)
By $\cO(\mbox{quantum space})$ we denote the polynomial unital $*$-algebra of a
quantum space defined by generators and relations,
 and by $C(\mbox{quantum space})$ the corresponding $C^*$-algebra. 
 
 The latter means the universal enveloping $C^*$-algebra (or the $C^*$-completion, or 
the $C^*$-closure) of a $*$-algebra
 in the following sense. Let $\cO$ be a unital $*$-algebra such that the set
Rep$(\cO)$ of all its bounded $*$-representations is non-empty, and such
that
\[
\forall\,x\in\cO:\,\sup\{ \|\varrho(x)\|\,|\, \varrho\in \mbox{Rep$(\cO)$}\}=:
\|x\|_{sup}<\infty.
\]
Then $\|\phantom{x}\|_{sup}$ defines a norm on the quotient algebra
$\cO/\{x\!\in\!\cO\,|\,\|x\|_{sup}=0\}$ whose norm completion is the aforementioned
universal enveloping $C^*$-algebra.

The thus constructed $C^*$-algebras indeed enjoy a very important universality
property which we use a number of times in this paper. One can show that,
if $C$ is the universal enveloping $C^*$-algebra of $\cO$ and $A$ is a
$C^*$-algebra,
 then for any $*$-homomorphism $\cO\ra A$ there exists a unique $C^*$-homomorphism 
$C\ra A$ such that the diagram 
 \[
\begin{diagram}[height=10mm]
\cO&\rTo&A\\
 &\rdTo&\uTo\\
 &&C
\end{diagram}
\]
 is commutative. (Here the diagonal map is the canonical quotient map combined with
 the injection of the quotient algebra into its norm completion.) This universality property
 is very useful because it allows one to define $C^*$-homomorphisms simply by appropriately
 specifying them on generators.
 
 On the other hand,  for $K$-theoretical computations it is useful to present a $C^*$-algebra
 as a fiber-product or crossed-product $C^*$-algebra. We denote the fiber product of 
 algebras by a decorated direct sum, and the crossed product of an algebra $A$ by a group
 $G$ as $A\rtimes G$. For generalities on the fiber products (pullbacks) of 
 $C^*$-algebras we recommend \cite{p-gk99}, 
 and for $C^*$-algebraic crossed products 
 we refer to \cite[Chapter V]{b-b98} and \cite[Chapter~7]{p-gk79}. Thanks to Theorem~\ref{fibreuni}, 
 our calculation of $K$-groups rests on the standard tools for  the $K$-theory of
 fiber products and crossed products.
%By classical points we understand 1-dimensional $*$-representations.

\subsection{A Heegaard splitting of \boldmath $S^3$}

Splitting and glueing topological spaces along 2-spheres or 2-tori are standard
procedures in the study of 3-dimensional manifolds. In the case of $S^3$, we have
the well-known Heegaard splittings. They present $S^3$
as two copies of a solid torus glued along their boundaries.

More precisely, the situation is as follows (see \cite[Section~0.3]{n-gl97} 
and \cite{hms} for related details).
Define 
\beq\label{class}
X=\{(z_1,z_2)\in\C^2\;|\;(1-|z_1|^2)(1-|z_2|^2)=0,\; |z_i|\leq 1\}.
\eeq
It is clear that $X$ is a glueing of two solid tori along their boundaries.
On the other hand,
 $S^3=\{(c_1,c_2)\in\C^2\;|\; |c_1|^2+|c_2|^2=1\}$
is the usual round 3-sphere.
To see that they are indeed homeomorphic, note that the following formulas
give continuous and mutually inverse maps between these two spaces:
\[\label{f}
f((z_1,z_2))=(|z_1|^2+|z_2|^2)^{-\frac{1}{2}}(z_1,{z_2}),
\]
\[\label{g}
g((c_1,c_2))=\frac{\sqrt{2}(c_1,{c_2})}{\sqrt{1+|\;|c_1|^2-|c_2|^2\;|}}\; .
\] 

It is a basic corollary of Bott periodicity \cite{b-pf72} that for the $n$-sphere $S^n$,
the Chern character 
\beq
\mathop{ch}: K^*(S^n)\lra H^*(S^n;\Q)
\eeq
maps $K^*(S^n)$ isomorphically onto $H^*(S^n;\Z)\inc H^*(S^n;\Q)$.
In particular, 
\beq
K^j(S^n)=
\left\{\begin{array}{cc}\Z\oplus\Z&j=0\\[.3cm]
0&j=1
\end{array}\right.,
~~~\mbox{for $n$ even, and}~~~ 
K^j(S^n)=
\left\{\begin{array}{cc}\Z&j=0\\[.3cm]
\Z&j=1
\end{array}\right.
~~~\mbox{for $n$ odd}.
\eeq
In the special case of $S^3$, any complex vector bundle on $S^3$ is trivial.
(This is because any bundle over $S^3$ can be clutched over $S^2$ from trivial bundles on 3-dimensional balls and 
the second homotopy group of any compact connected Lie group is zero.)
Thus $K^0(S^3)=\Z$ and a generator is the trivial line bundle $S^3\times\C$.
For $K^1(S^3)$, observe that, by Bott periodicity, $\pi_3(U(n))=\Z$ for $n\geq 2$.
Hence $K^1(S^3)=\Z$ with  a generator given by the map
\beq
(z_1,z_2)\longmapsto\left(\ba{cc}
z_1&-\overline{z}_2\\
z_2&\overline{z}_1\\
\ea\right).
\eeq

From the $C^*$-algebra point of view, $K_0(C(S^3))$ is generated by the unit of
$C(S^3)$ and $K_1(C(S^3))$ is generated by the 2-by-2 unitary matrix given
above, now with $z_1$, $z_2$, viewed as the coordinate functions on $S^3$.
If we view $S^3$ as the union of two solid tori, glued along their boundaries,
then there is the 6-term Mayer-Vietoris exact sequence 
\beq\label{classmv}
\begin{diagram}[height=10mm]
K_0(C(S^3))  &  \rTo 
& K_0(C(D^2\!\times\! S^1))\oplus K_0(C(S^1\!\times\! D^2)) &\rTo  
&  K_0(C(S^1\!\times\! S^1)) \\
 \uTo &  &  & &\dTo \\
K_1(C(S^1\!\times\! S^1))& \lTo     
&K_1(C(D^2\!\times\! S^1))\oplus K_1(C(S^1\!\times\! D^2)) &\lTo 
&K_1(C(S^3))\; . \\ 
\end{diagram}
\eeq
The $K$-groups of $S^3$ could be calculated from this diagram. In the noncommutative
case, this is how we shall calculate the $K$-groups for our examples of quantum 3-spheres.

\subsection{The Toeplitz and irrational rotation \boldmath$C^*$-algebras}

The Toeplitz and irrational rotation $C^*$-algebras ${\mathcal T}$ and 
$C(S^1) {\rtimes}_ {\theta}\,\Z$ are among the best studied examples of operator algebras.
They admit natural geometrical interpretations as a quantum disc and a quantum torus,
respectively. As explained in the introduction, they become basic building blocks in our construction of a quantum 3-sphere
just as a disc and torus can be used to build up~$S^3$. Herein, we recall some basic
facts about these algebras.

Let $\cO(D_q)$ stand for the unital $*$-algebra generated by $z,z^*$ satisfying
\beq\label{d}
z^*z-qzz^*=1-q,~~~0\leq q<1.
\eeq
Denote by $\{e_k\}_{k\in\N}$ an orthonormal basis of a Hilbert space.
Up to the unitary equivalence,  $\cO(D_q)$ has 
the following irreducible $*$-representations in
bounded operators \cite{kl93}:
\beq\label{r}
\pi(z)\,e_k=\sqrt{1-q^{k+1}}\;e_{k+1},\; k\in\N,~~~\pi_\lambda(z)=\lambda,\;\lambda\in U(1).
\eeq
Using this classification, one can show that the extension of the representation $\pi$  to the universal enveloping
$C^*$-algebra $C(D_q)$ of $\cO(D_q)$ is faithful \cite[p.14]{kl93}.
A straightforward argument  proves  that  the norm of $z$ in  $C(D_q)$  is~1.
 (The same argument is used
to prove (\ref{norm}).)
For $q=1$ this norm condition is no longer a consequence of the commutation relation,
 but can be naturally
assumed to hold, as
it always holds automatically for $0\leq q<1$. The norm condition combined with the commutativity of $z$
and $z^*$ yields the $C^*$-algebra of continuous functions on the unit disc~$D$.
This is a motivation for calling $C(D_q)$ the $C^*$-algebra of a quantum disc~\cite{kl93}.

On the other hand, $C(D_0)$ is, by construction, the well-known Toeplitz algebra~$\cT$.
The faithful representation $\pi$ allows us to view $\cT$ as the operator algebra generated by
the unilateral shift.
It appears quite interesting that, for $0<q<1$, we can always re-scale the generator $z$
by dividing it (on the right) by its absolute value to obtain the usual generator of $\cT$
that we have for $q=0$. This is precisely the phenomenon that we use to prove
Theorem~\ref{00pq} (see (\ref{gst})).

Next, let $\cO(T^2_\theta)$ denote the unital $*$-algebra generated by $x$ and $y$
satisfying
\[\label{t2t}
xx^*=1=x^*x,~~~yy^*=1=y^*y,~~~xy=e^{2\pi i\theta}yx,~~~\theta\in [0,1[\,.
\]
The universal enveloping $C^*$-algebra $C(T^2_\theta)$ of $\cO(T^2_\theta)$ 
coincides with the crossed-product algebra \mbox{$C(S^1)\rtimes_\theta\,\Z$}. 
For $\theta=0$ we get the algebra of continuous functions on the 2-torus~$T^2$. Hence,
for $\theta\neq0$, the quantum space $T^2_\theta$ is called a noncommutative torus.
The situation is particularly interesting when $\theta$ is irrational. Therefore, the irrationality
of $\theta$ is often automatically assumed, in which case the crossed product
$C(S^1)\rtimes_\theta\,\Z$ goes under the name irrational rotation $C^*$-algebra 
\cite{r-ma90}.

\subsection{The Mayer-Vietoris sequence for the K-theory of C*-algebras}

Let
\beq
\begin{CD}
A @>>{}> B_2 \\
@VV{}V @VV{\pi_2}V \\
B_1 @>>{\pi_1}> D
\end{CD}
\eeq
be a commutative diagram of $C^*$-algebras in which each 
homomorphism is surjective and $A := \{ (b_1 , b_2) \in B_1 \oplus B_2 \mid \pi_1(b_1) = \pi_2(b_2) \}$. Then $A$ is
again a $C^*$-algebra and it is called the fiber product (or pullback) of $B_1$ and $B_2$ 
over $D$. 
%(We refer to \cite{p-gk99} for generalities on this construction.)
In case the $C^*$-algebras $B_1$ and $B_2$ are commutative, 
an argument of Atiyah and Hirzebruch \cite[p.32]{ah62}, allows one to
construct the six-term $K$-theory 
exact sequence
\beq\label{mv}
\begin{CD}
K_0 (A) @>>{}> K_0 (B_1) \oplus K_0 (B_2) @>>{}> K_0 (D) \\
@AA{}A &&@VV{}V \\
K_1 (D) @<<{}< K_1 (B_1) \oplus K_1 (B_2) @<<{}< K_1 (A)\, .
\end{CD}
\eeq
We make a very slight change in the Atiyah-Hirzebruch argument 
(so that no partition of unity is used) and then observe that their construction 
applies to the case when the $C^*$-algebras are not required to be commutative.

First, recall the classical  construction. Here $X$ is a locally compact Hausdorff
 space and $X_1$, $X_2$, are closed subsets of $X$ with $X = X_1 \cup X_2$. Let $\tilde X$ be the subset of $X \times [0,1]$:
\[
\tilde X = 
X_1 \!\times\! \{ 0 \} \;\cup\; (X_1 \cap X_2)\! \times\! [0,1] \;\cup\; X_2 \!\times\! \{ 1 \} \, .
\]
Thus $\tilde X$ is the disjoint union of $X_1$ and $X_2$ with a cylinder connecting the two copies of $X_1 \cap X_2$
\vspace*{-5mm}\begin{figure}[h]
$$
\includegraphics[width=50mm]{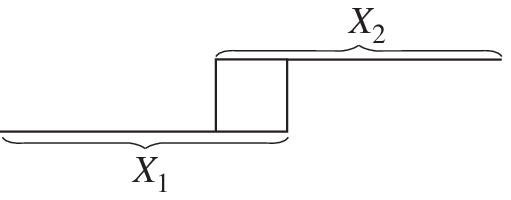}
$$
\end{figure}\vspace*{-7mm}

\noindent
that is topologized as a subspace of $X \times [0,1]$. Consider now the six-term $K$-theory 
exact sequence
 \cite[p.90]{k-m78}
 for $(\,\tilde X\, ,\, X_1 \!\times\! \{ 0 \}\: \cup\: X_2 \!\times\! \{ 1 \}\,)$:
\beq\label{mfatiyah}
\begin{diagram}[height=10mm]
K^0 (\tilde X\, ,\, X_1 \!\times\! \{ 0 \} \:\cup\: X_2 \!\times\! \{ 1 \})  &  \rTo 
& K^0 (\tilde X)  &\rTo  
&  K^0 (X_1) \oplus K^0 (X_2)  \\
 \uTo &  &  & &\dTo \\
K^1 (X_1) \oplus K^1 (X_2) & \lTo     
&  K^1 (\tilde X) &\lTo 
& K^1 (\tilde X\, , \,X_1\! \times\! \{ 0 \}\: \cup\: X_2 \!\times\! \{ 1 \})\; . \\ 
\end{diagram}
\eeq
Next, let us use the notation $Y/Y'$ to signify the space $Y$ with its subspace $Y'$ shrunk
to a point and observe that $\tilde X\, / \,(X_1 \!\times\! \{ 0 \} \cup X_2 \!\times\! \{ 1 \})$ 
is the suspension of $X_1 \cap X_2$. Hence we conclude that
\beq\label{susp}
K^j (\tilde X\, ,\, X_1\! \times\! \{ 0 \} \cup X_2 \!\times\! \{ 1 \}) = K^{1-j} (X_1 \cap X_2) \, .
\eeq
Furthermore, we can prove that
  the restriction of the projection  $X \times [0,1] \to X$ to $\tilde X \to X$  gives an isomorphism in $K$-theory: 
\[\label{iso}
K^j (X) \overset{\cong}{\longrightarrow} K^j (\tilde X),~~~  j=0,1.
\]
Indeed, the map $\tilde X \to X$ maps $X_1 \!\times\! \{ 0 \} \;\cup\; (X_1 \cap X_2) \!\times\! [0,1]$ to $X_1$ and so gives a map of pairs
\[
(\tilde X \,,\, X_1\! \times\! \{ 0 \} \,\cup\, (X_1 \cap X_2) \!\times\! [0,1]) \to (X,X_1) \, .
\]
Hence the six-term $K$-theory exact sequence for $(X,X_1)$ maps to the six-term $K$-theory exact sequence for $(\tilde X\, ,\, X_1 \!\times\! \{ 0 \}\, \cup\, (X_1 \cap X_2)\! \times \![0,1])$. Since $\tilde X\, /\, (X_1 \!\times\! \{ 0 \} \,\cup\, (X_1 \cap X_2)\! \times\! [0,1]) 
= X/X_1$ and since $X_1 \!\times\! \{ 0 \} \,\cup\, (X_1\cap X_2)\! \times\! [0,1] \to X_1$ is clearly a homotopy equivalence, the Five Isomorphisms Lemma now applies to end the
proof of (\ref{iso}). Combining this with (\ref{susp}) and the six-term exact sequence for 
$(\tilde X\, ,\, X_1 \!\times\! \{ 0 \} \:\cup\: X_2\! \times\! \{ 1 \})$ now yields the desired Mayer-Vietoris sequence (\ref{mv}). 

 We shall now show 
  that the reasoning used in the commutative case applies basically unchanged to the general case.
Let $\tilde A \inc B_1 \oplus B_2 \oplus C ([0,1],D)$ be
\[
\tilde A = \{ (b_1 , b_2 , \omega) \mid b_1 \!\in\! B_1 ,\, b_2 \!\in\! B_2 ,\,
 \omega (0) = \pi_1 (b_1) ,\, \omega (1) = \pi_2 (b_2)\} \, .
\]
As usual, $C ([0,1] , D)$ is the $C^*$-algebra of all continuous maps from $[0,1]$ to $D$.
Similarly, put
\[
C_0 (\,]0,1[\,,D) = \{ \omega \in C ([0,1] , D) \mid \omega (0) = \omega (1) = 0 \} \, .
\]
Consider now the short exact sequence of $C^*$-algebras
\[
0 \lra C_0 (\,]0,1[\,,D) \lra \tilde A \lra B_1 \oplus B_2 \lra 0,\;\;\;
(b_1 , b_2 , \omega) \longmapsto (b_1 , b_2) \, .
\]
The six-term $K$-theory exact sequence of this short exact sequence of $C^*$-algebras plus the identifications
\[
K_j (C (\,]0,1[\,,D)) = K_{1-j} (D),\;\;\; K_j (\tilde A )= K_j (A), \;\;\; j=0,1,
\]
then yield the desired Mayer-Vietoris exact sequence. The first of these identifications is immediate because $C(\,]0,1[\,,D)$ is the suspension of $D$. For the second, we wish to prove that the evident map
\[
A \lra \tilde A,\;\;\;
(b_1 , b_2) \longmapsto (b_1 , b_2 , \underline{\pi_1 (b_1)}),
\]
where $\underline{\pi_1 b_1}$ is the constant path at $\pi_1 (b_1) = \pi_2 (b_2)$,
 gives a $K$-theory isomorphism.
Let $I_1$ be the kernel of
\[
A \lra B_1,\;\;\;
(b_1 , b_2) \longmapsto b_1 \, .
\]
Consider the pair $(A,I_1)$, and note that $I_1 \cong\ker(\pi_2 : B_2 \to D)$. 
The image of $I_1$ in $\tilde A$ is $\{ (0,b_2 , \underline{0}) \mid \pi_2 (b_2) = 0 \}$. (Here $\underline{0}$ denotes the constant path at $0 \in D$.) Thus the image of $I_1$ in $\tilde A$ is an ideal in $\tilde A$. By a slight abuse of notation, denote the image of $I_1$ in $\tilde A$ also by $I_1$. The map of pairs 
\[
(A,I_1) \lra (\tilde A , I_1),\;\;\; (b_1 , b_2) \longmapsto (b_1 , b_2 , \underline{\pi_1 (b_1)}),
\]
maps the six-term $K$-theory exact sequence
\cite[p.67]{b-b98}
 for $(A,I_1)$ to the six-term $K$-theory exact sequence for $(\tilde A , I_1)$.
Finally, let us compare $A/I_1 = B_1$ and $\tilde A / I_1$. We have
\begin{eqnarray}
\tilde A / I_1 &= &\{ (b_1 , b_2 , \omega) \mid \omega (0) = \pi_1 (b_1) , \;\omega (1) =
 \pi_2 (b_2) \}\; /\; \{ (0,b_2 , \underline{0}) \mid \pi_2 (b_2) = 0 \} \nonumber \\
&= &\{ (b_1 , \omega) \mid b_1 \in B_1 \ \hbox{and} \ \omega : [0,1] \to D \ \hbox{has} \ \omega (0) = \pi_1 (b_1) \} \, . 
\end{eqnarray}
It is now clear that the $C^*$-algebras $A/I_1 = B_1$ and $\tilde A / I_1$ are homotopy equivalent. Therefore, an analogue of the argument proving (\ref{iso})
 applies to prove the desired identification $K_j (\tilde A )\cong K_j ( A)$,  $j=0,1$.
 This completes the construction of the Mayer-Vietoris exact sequence for
 the $K$-theory of $C^*$-algebras.

\section{The universal C*-algebras}
% of the Heegaard quantum 3-spheres}
\setcounter{equation}{0}

In this section, we show the existence of the universal  $C^*$-algebra for relations that
we interpret as defining a 3-parameter family of noncommutative 3-spheres. By definition,
the topology of the parameter space is that of a solid torus, but it turns out that
 identifying points in this space that give isomorphic $C^*$-algebras largely reduces its size.
We begin by defining the polynomial $*$-algebras of our noncommutative 3-spheres
\cite{bhms}.
\bde
 The  unital $*$-algebra $\cO(S^3_{pq\theta})$ is the quotient 
of the free $*$-algebra generated by elements $a$ and $b$  by the
$*$-ideal generated by the following relations:
\begin{align}
\label{sphere}
&(1-aa^*)(1-bb^*)=0 & &\mbox{the sphere equation},
\\ \label{direl}
&a^*a=paa^*+1-p,~~~b^*b=qbb^*+1-q & &\mbox{quantum-disc relations},
\\ \label{torel}
&ab=e^{2\pi i\theta} ba, ~~~ab^*=e^{-2\pi i\theta}b^*a
& & \mbox{noncommutative-torus relations}.
\end{align}
Here $p,q$ and $\theta$ are parameters running the unit interval $[0,1]$.
\ede\noindent
It turns out to be very helpful to distinguish the following elements of this algebra:
\[\label{AB}
 A:=1-aa^*\;\;\;\mbox{and}\;\;\;B:=1-bb^*.
\]
They satisfy $AB=0$ and one can think of them as some key
functions supported, respectively, on one or the other noncommutative solid torus whose glueing yields our quantum 3-sphere. It is also immediate that we have the following 
relations:
\beq
Ab=bA,~~~Ba=aB,~~~Aa=paA,~~~Bb=qbB,~~~a^*a=1-pA,~~~b^*b=1-qB.
\eeq

To prove the existence  of the enveloping $C^*$-algebra
of  $\cO(S^3_{pq\theta})$, we first need to study bounded $*$-representations of the latter.
It is straightforward to verify that for all values of parameters $p$, $q$ and $\theta$ we have:
\bpr\label{reps}
Let $\{e_k\}_{k\in\N}$ and $\{e_{m,n}\}_{m\in\Z,n\in\N}$ be
orthonormal bases of two
Hilbert spaces, respectively. Put $\mu:=e^{2\pi i\theta}$. Then  the
following formulas define bounded $*$-representations
of~$\cO(S^3_{pq\theta})$:
\vspace*{-5mm}
\begin{itemize}\item[]\begin{itemize}\item[]\begin{itemize}
\vspace*{4mm}\item[(1)]~
%a family $\rho_{\lambda}$, $\hl\in U(1)$, given by
%\beq
$
\rho_{\lambda}(a)\,e_{k}=\sqrt{1-p^{k+1}}\;e_{k+1},
~~~\rho_{\lambda}(b)\,e_{k}=\hl\mu^{-k}\;e_{k},~~~\hl\in U(1),
%~~~\mu:=e^{2\pi i\theta},
%\eeq
$
\vspace*{4mm}\item[(2)]~
%a family $\rho'_{\lambda}$, $\hl\in U(1)$, given by
%\beq
$
\rho'_{\lambda}(a)\,e_{k}=\hl\mu^{k}\;e_{k},
~~~\rho'_{\lambda}(b)\,e_{k}=\sqrt{1-q^{k+1}}\;e_{k+1},~~~\hl\in U(1),
%~~~\mu:=e^{2\pi i\theta},
%\eeq
$
\vspace*{4mm}\item[(3)]~
% $\rho$ given by
%\beq
$
\rho(a)\,e_{m,n}=\mu^m\sqrt{1-p^{n+1}}\;e_{m,n+1},~~~
\rho(b)\,e_{m,n}=e_{m+1,n},
%~~~\mu:=e^{2\pi i\theta},
%\eeq
$
\vspace*{4mm}\item[(4)]~
% $\rho'$ given by
%\beq
$
\rho'(a)\,e_{m,n}=e_{m+1,n},~~~
\rho'(b)\,e_{m,n}=\mu^{-m}\sqrt{1-q^{n+1}}\;e_{m,n+1}.
%~~~\mu:=e^{2\pi i\theta}.
%\eeq
$
\end{itemize}\end{itemize}\end{itemize}
\epr\bigskip

\noindent
In order to check the relations, we need explicit formulas for the adjoints of the 
aforelisted operators. They are as follows.
\vspace*{3mm}
\begin{align}
%\nonumber\\
\rho_{\lambda}(a^*)\,e_{k}=
\left\{\ba{lcc} 
\!\!\!\!\!\!& \sqrt{1-p^k}\;e_{k-1} & \mbox{ for } k>0\\
\!\!\!\!\!\!& 0 & \mbox{ for } k=0
\ea\right.,
~~~~~~~
\rho_{\lambda}(b^*)\,e_{k}=\hl^{-1}\mu^k\;e_{k}\,,
\\  \nonumber\\ 
\rho'_{\lambda}(a^*)\,e_{k}=\hl^{-1}\mu^{-k}\;e_{k}\,,
~~~~~~
\rho'_{\lambda}(b^*)\,e_{k}=
\left\{\ba{lcc} 
\!\!\!\!\!\!& \sqrt{1-q^k}\;e_{k-1} & \mbox{ for } k>0\\
\!\!\!\!\!\!& 0 & \mbox{ for } k=0
\ea\right.,
\\  \nonumber\\
\rho(a^*)\,e_{m,n}=
\left\{\ba{lcc} 
\!\!\!\!\!\!& \mu^{-m}\sqrt{1-p^n}\;e_{m,n-1} & \mbox{ for } n>0\\
\!\!\!\!\!\!& 0 & \mbox{ for } n=0
\ea\right.,
~~~~~~
\rho(b^*)\,e_{m,n}=e_{m-1,n}\,,
\\ \nonumber \\ 
\rho'(a^*)\,e_{m,n}=e_{m-1,n}\,,
~~~~~~
\rho'(b^*)\,e_{m,n}=
\left\{\ba{lcc} 
\!\!\!\!\!\!& \mu^{m}\sqrt{1-q^n}\;e_{m,n-1} & \mbox{ for } n>0\\
\!\!\!\!\!\!& 0 & \mbox{ for } n=0
\ea\right..
%\\  \nonumber
\end{align}~\\
With these formulas at hand, one immediately computes:
\vspace*{3mm}
\begin{align}
%\nonumber\\
&
\rho_{\lambda}(A)\,e_{k}=
\left\{\ba{lcc} 
\!\!\!\!\!\!&p^k\;e_{k} & \mbox{ for } k>0\\
\!\!\!\!\!\!& e_{k} & \mbox{ for } k=0
\ea\right.,
~~~~~~~
\rho_{\lambda}(B)=0,
\\ \nonumber\\ &
\rho'_{\lambda}(A)=0, 
~~~~~~~
\rho'_{\lambda}(B)\,e_{k}=
\left\{\ba{lcc} 
\!\!\!\!\!\!&q^k\;e_{k} & \mbox{ for } k>0\\
\!\!\!\!\!\!& e_{k} & \mbox{ for } k=0
\ea\right.,
\\ \nonumber\\
&
\rho(A)\,e_{m,n}=
\left\{\ba{lcc} 
\!\!\!\!\!\!&p^n\;e_{m,n} & \mbox{ for } n>0\\
\!\!\!\!\!\!& e_{m,n} & \mbox{ for } n=0
\ea\right.,
~~~~~~
\rho(B)=0,
\\ \nonumber\\&
\rho'(A)=0,
~~~~~~
\rho'(B)\,e_{m,n}=
\left\{\ba{lcc} 
\!\!\!\!\!\!&q^n\;e_{m,n}& \mbox{ for } n>0\\
\!\!\!\!\!\!& e_{m,n} & \mbox{ for } n=0
\ea\right..
\\ \nonumber
\end{align}
\note{
Note that there are no classical points of $\cO(S^3_{pq\theta})$:
It follows from (\ref{direl}) that $a$ and $b$ would be represented
by numbers of modulus one, contradicting the fact that due to
(\ref{torel}) one of these numbers would have to vanish.
Notation: $a^k:={a^*}^{-k}$ for $k<0$.
}
\indent%\hspace*{5mm}
Note now that, for $0\leq p,q<1$, the relations (\ref{direl}) entail that the norm of $a$ and $b$
in any bounded $*$-representation $\varrho$ is at most 1. Indeed, we have
\[\label{1-p}
1-p = \|\varrho(a^*a-paa^*)\| \geq\|\varrho(a^*a)\|-p\|\varrho(aa^*)\| = (1-p)\|\varrho(a)\|^2,
\]
which implies $\|\varrho(a)\|\leq 1$. The same reasoning works for $\|\varrho(b)\|$.
Hence, as  the space of bounded $*$-representations of $\cO(S^3_{pq\theta})$ is non-empty
and the generators are uniformly bounded in all of them, there exists the universal
enveloping $C^*$-algebra of $\cO(S^3_{pq\theta})$.
\bde\label{cdef}
For $0\leq p,q<1$, the $C^*$-algebra of $S^3_{pq\theta}$ is the universal enveloping 
$C^*$-algebra of $\cO(S^3_{pq\theta})$. We denote it by $C(S^3_{pq\theta})$.
\ede 
\noindent
To finish the discussion of the norm of the generators, observe that 
the formulas Proposition~\ref{reps}(3) entail $\|\rho(a)\|,\|\rho(b)\|\geq 1$.
This fact combined with
the above shown inequalities $\|\varrho(a)\|,\|\varrho(b)\|\leq 1$, for any $\varrho$,
implies that
\[\label{norm}
\|a\|=1=\|b\|,~~~\forall\; 0\leq p,q<1.
\]

We skip the analysis of the irreducibility of the above representations and postpone
the discussion of what happens for special values of parameters $p$, $q$ and $\theta$
for later on. Instead, we focus on the representations $\rho$ and $\rho'$ and find a linear basis of $\cO(S^3_{pq\theta})$ for two cases of parameter values.
Contrary to representation formulas, basis expressions 
 are sensitive to special parameter values.
To write basis elements, it is convenient to adopt the following convention:
\[\label{conv}
x_\alpha=
\left\{\ba{lccc} 
\!\!\!\!\!\!& x^\alpha& \mbox{ for } &\alpha\in\N\\
\!\!\!\!\!\!& (x^*)^{-\alpha} & \mbox{ for } &\alpha\in\Z\setminus\N
\ea\right..
\]
\ble\label{baspq}
 If $0<p,q<1$,
then the elements $a_\alpha b_\beta A^m$, $a_{\alpha'} b_{\beta'}B^{m'}$, $\alpha,\ha',\hb,\hb',m,m'\in\Z$, $m\geq0$, $m'>0$,
form a linear basis of $\cO(S^3_{pq\theta})$, and 
the representation $\rho\oplus\rho'$ is faithful.
\ele
\bpf
It is a straightforward consequence of  relations in the algebra
$\cO(S^3_{pq\theta})$ that the 
aforelisted elements span this algebra.
Now, let
\beq
f:=\sum_{\alpha,\beta \in\Z,\,l\in\N} f_{\alpha \beta l}\,a_\alpha b_\beta A^l
\;+\sum_{\alpha ,\beta \in\Z,\,l>0} g_{\alpha \beta l}\,a_\alpha b_\beta B^l
\eeq
be an arbitrary linear combination, and assume that it is equal to~0.
Then $\rho(f)e_{m,n}=0$ for any $m\in\Z$, $n\in\N$. On the other hand,

\beq\label{bro2}
\rho(f)\,e_{m,n}=\!\!\!\sum_{\alpha ,\beta \in\Z,\,l\in\N}\!\!\!
f_{\alpha \beta l}\,p^{nl}\mu^{m\alpha} \left\{\ba{llllccc}
\!\!\!&\sqrt{(1-p^{n+1})\cdots(1-p^{n+\alpha })}\;e_{m+\beta ,n+\alpha }
&\mbox{for}&
\alpha > 0\\
\!\!\!\!\!\!&&&\\
\!\!\!\!\!\!&e_{m+\beta ,n+\alpha }&\mbox{for}&\alpha =0\\
\!\!\!\!\!\!&&&\\
\!\!\!\!\!\!&\sqrt{(1-p^n)\cdots(1-p^{n+\alpha +1})}\;e_{m+\beta ,n+\alpha }
&\mbox{for}&
\alpha <0
\ea\right..
\eeq

\noindent
Since $f$ is a finite linear combination, there exists 
$\alpha _0=\min\{\alpha\in\Z\, |\,f_{\alpha \beta l}\neq0,
\beta\in\Z,\,l\in\N\}$. Hence,  the square roots in (\ref{bro2})
are non-zero for all $n\geq -\alpha _0$. Therefore, it follows from the linear independence
of the basis vectors that
\beq
\sum_{l\in\N}f_{\alpha \beta l}(p^n)^l=0,~~~\forall\, \alpha ,\beta \in\Z,\;n\geq -\alpha _0\,.
\eeq
Thus $f_{\alpha \beta l}$
are coefficients of a polynomial vanishing at infinitely many points. 
Consequently $f_{\alpha \beta l}=0$ for all $\alpha ,\beta ,l$.

Applying now $\rho'(f)$ to $e_{m,n}$ we obtain

\beq%\label{bro2}
\rho'(f)\,e_{m,n}=\!\!\!\sum_{\alpha ,\beta \in\Z,\,l\in\N}\!\!\!
g_{\alpha \beta l}\,q^{nl}\mu^{-m\beta} \left\{\ba{llllccc}
\!\!\!&\sqrt{(1-q^{n+1})\cdots(1-q^{n+\beta })}\;e_{m+\alpha ,n+\beta }
&\mbox{for}&
\beta > 0\\
\!\!\!\!\!\!&&&\\
\!\!\!\!\!\!&e_{m+\alpha ,n+\beta }&\mbox{for}&\beta =0\\
\!\!\!\!\!\!&&&\\
\!\!\!\!\!\!&\sqrt{(1-q^n)\cdots(1-q^{n+\beta +1})}\;e_{m+\alpha ,n+\beta }
&\mbox{for}&
\beta <0
\ea\right..
\eeq\\

\noindent
Using the same argument as above, we can conclude that, for all sufficiently big $n\in\N$
and for all $\ha,\hb\in\Z$,
\beq
\sum_{l>0}g_{\alpha \beta l}(q^n)^l=0.
\eeq
Again as above, this entails the vanishing of all $g_{\alpha \beta l}$.

Finally, let us note that we have shown that $(\rho\oplus\rho')(f)=0$ implies
$f=0$. This means the desired faithfulness of $\rho\oplus\rho'$.
\epf

\bre
Let $z$ be the generator of $\cO(D_q)$ (see (\ref{d})). Assume the convention (\ref{conv}), and put $Z:=1-zz^*$. 
For $0<q<1$, reasoning much as in the proof of the above lemma, one can show that
the 
elements $z_\alpha Z^k$, $\ha\in\Z$, $k\in\N$, form a basis of
$\cO(D_q)$.
\ere

Let us now pass to considering the case $p=0=q$. Our main point here is that,
for $0\leq p,q<1$, we have $C(S^3_{pq\theta})\cong C(S^3_{00\theta})$. To avoid confusion,
we shall denote the generators of $\cO(S^3_{00\theta})$ by $s$ and $t$, respectively,
rather than by $a$ and $b$. For the time being, the latter notation is reserved for the cases
$0<p,q<1$. For $p=0=q$ the relations (\ref{direl}) become simpler and the elements 
$A$ and $B$ (\ref{AB}) become projections. Hence, it is clear that their powers can
no longer be useful in writing down basis elements. Therefore, we introduce the following
 self-adjoint elements:
\beq\label{akbk}
A_k=1-s^k{s^*}^k,~~~B_k=1-t^k{t^*}^k,~~~k\in\N.
\eeq
For these elements, one can immediately derive the identities:
\beq\label{at}
[t,A_k]=0,~~~A_{k+1}s=sA_{k},~~~s^*A_{k+1}=A_{k}s^*,~~~A_{k+1}=sA_ks^*+A_1\,,
\eeq
\beq\label{bs}
[s,B_k]=0,~~~B_{k+1}t=tB_{k},~~~t^*B_{k+1}=B_{k}t^*,~~~B_{k+1}=tB_{k}t^*+B_1\,.
\eeq
With this notation, we can write the sphere relation (\ref{sphere}) as $A_1B_1=0$. 
Now, combining the first formula of (\ref{bs}) and the last of (\ref{at}) with the standard 
induction proves that $A_kB_1=0$, $\forall\,k\in\N$. Applying the standard 
induction with this formula as the starting point, and the first formula of (\ref{at}) and the
last of (\ref{bs}) for the induction step, yields
\[\label{akbl0}
A_kB_l=0,~~~\forall\,k,l\in\N.
\]

Again adopting the convention (\ref{conv}), we can now claim:
\ble\label{bas3}
Let $\ha,\beta,k\in\Z$. The elements $s_\alpha t_\beta$, and  $s_\alpha t_\beta A_k$ with $k>0$, $\alpha>-k$,
and $s_\alpha t_\beta B_k$ with $k>0$, $\beta>-k$,  form a linear basis of $\cO(S^3_{00\theta})$. The representation $\rho\oplus\rho'$ is faithful.
\ele
\bpf
It follows directly from relations in $\cO(S^3_{00\theta})$ that the monomials
 $s_\mu t^k{t^*}^l$, $\mu\in\Z$, $k,l\in\N$,
and $s^m{s^*}^nt_\nu$, $\nu\in\Z$, $m,n\in\N$, span $\cO(S^3_{00\theta})$.
We now want to show that so do the elements listed in the assertion of the lemma.
Using (\ref{bs}), for $l>k>0$, we obtain
\beq
s_\mu t^k{t^*}^l=s_\mu t^k{t^*}^k{t^*}^{l-k}=s_\mu(1-B_k){t^*}^{l-k}=s_\mu {t^*}^{l-k}-
s_\mu {t^*}^{l-k}B_l\,.
\eeq 
Analogously, for $k\geq l$, we have
\beq
s_\mu t^k{t^*}^l=s_\mu t^{k-l}t^l{t^*}^l=s_\mu t^{k-l}(1-B_l)=s_\mu t^{k-l}-s_\mu t^{k-l}B_l\,.
\eeq
Similarly, taking an advantage of (\ref{at}) one checks easily for $m\geq n$ that
\beq
s^m{s^*}^nt_\nu=s^{m-n}s^n{s^*}^nt_\nu=s^{m-n}t_\nu-s^{m-n}t_\nu A_n\,.
\eeq
For $n>m>0$, it follows from (\ref{at}) that
\beq
s^m{s^*}^nt_\nu=s^m{s^*}^m{s^*}^{n-m}t_\nu={s^*}^{n-m}t_\nu-{s^*}^{n-m}t_\nu A_n\,.
\eeq
Summarizing, we have shown that the elements given in the assertion of the lemma
span~$\cO(S^3_{00\theta})$.

Let us now write a general element $f\in\cO(S^3_{00\theta})$ as a linear combination:
%of these elements.
\beq
f:=\sum_{\stackrel{\mbox{\scriptsize$\alpha,\beta,k\in\Z$}}{\alpha>-k<0}}
f_{\alpha\beta k}\,s_\alpha t_\beta A_k\;
+\sum_{\stackrel{\mbox{\scriptsize$\alpha,\beta,k\in\Z$}}{\alpha>-k<0}}
g_{\alpha\beta k}\,s_\alpha t_\beta B_k\;
+\sum_{\alpha,\beta\in\Z}h_{\alpha,\beta}\,s_\alpha t_\beta.
\eeq
Using the straightforward to verify formulas
\begin{align}
\label{Ai}&
\rho(A_k)\,e_{m,n}=\left\{\ba{cc} e_{m,n}&\mbox{ for } k>n\\0&\mbox{ for } k\leq n\ea\right.,
~~~~~~\rho(B_k)\,e_{m,n}=0,\\
\label{Bi}&
\rho'(A_k)\,e_{m,n}=0,~~~~~~
\rho'(B_k)\,e_{m,n}=\left\{\ba{cc} e_{m,n}&\mbox{ for } k>n\\0&\mbox{ for } k\leq n\ea\right.,
\end{align}
we can compute the operators $\rho(f)$ and $\rho'(f)$.
For any $m\in\Z$ and $n\in\N$, employing the shorthand notation $\mu:=e^{2\pi i\theta}$,
we have
\begin{align}
&
\rho\;\Lblp\!\!\!\sum_{\stackrel{\mbox{\scriptsize$\alpha,\beta,k\in\Z$}}{\alpha>-k<0}}
f_{\alpha\beta k}\,s_\alpha t_\beta A_k\;
+\sum_{\stackrel{\mbox{\scriptsize$\alpha,\beta,k\in\Z$}}{\beta>-k<0}}
g_{\alpha\beta k}\,s_\alpha t_\beta B_k\;
+\sum_{\alpha,\beta\in\Z}h_{\alpha,\beta}\,s_\alpha t_\beta
\Lbrp e_{m,n}
\nonumber\\ &
=\sum_{\stackrel{\mbox{\scriptsize$\alpha,\beta,k\in\Z$}}
{-k<-n\leq\alpha}}
f_{\alpha\beta k}\,\mu^{(m+\beta)\alpha}\,e_{m+\beta,n+\alpha}\;
+\sum_{\stackrel{\mbox{\scriptsize$\alpha,\beta\in\Z$}}{\alpha\geq-n}}
h_{\alpha\beta}\,\mu^{(m+\beta)\alpha}\,
e_{m+\beta,n+\alpha}\,.\label{fh}
\end{align}

To prove the linear independence of our linear generators, we assume now that $f=0$.
We begin by drawing conclusions from the entailed equality $\rho(f)=0$.
If the first sum in (\ref{fh}) is not vanishing form the start, then there exists 
\[
k_0:=\max\{k\in\Z\,|\,f_{\alpha\beta k}\neq 0,\,\alpha,\beta\in\Z,\,\alpha>-k<0\},
\]
 and  choosing $n\geq k_0$ makes it zero.
Thus the first sum is zero at least for all sufficiently big~$n$, and  we obtain
\beq
\sum_{\stackrel{\mbox{\scriptsize$\alpha,\beta\in\Z$}}{\alpha\geq-n}}
h_{\alpha\beta}\,\mu^{(m+\beta)\alpha}\,
e_{m+\beta,n+\alpha}=0.%,\;\Lall m,n\in\Z,\;n\geq k_0\,.
\eeq
%If $n$ is also bigger than $\max\{-\ha\,|\,h_{\alpha\beta}\neq 0\}$
Again for all sufficiently big~$n$, the condition
$\ha\geq -n$ in the foregoing sum is automatically satisfied for all possibly non-vanishing coefficients $h_{\alpha\beta}$ that appear therein. Thus we can conclude
from the linear independence of the vectors $e_{m+\beta,n+\alpha}$ that
\beq
h_{\alpha\beta}=0, ~~~\forall\;\alpha,\beta\in\Z.
\eeq
Putting this back into (\ref{fh}) and noticing that terms with different $\beta$ are linearly independent, we get
\beq\label{fsum}
\sum_{\stackrel{\mbox{\scriptsize$\alpha,k\in\Z$}}
{-k<-n\leq\alpha}}
f_{\alpha\beta k}\,\mu^{(m+\beta)\alpha}\,e_{m+\beta,n+\alpha}=0,
~~~\Lall \hb,m\in\Z,n\in\N.
\eeq
Suppose now that at least one of the coefficients $f_{\alpha\beta k}$ is non-zero. Then
$k_0$ exists. Due to the assumption that $k$ is always a positive natural number,  
$k_0-1\geq 0$, so that we can
choose $n=k_0-1$. Then
the condition $k>k_0-1$ in (\ref{fsum}) entails that $k_0$ is the only possible value for $k$,
and that there is no summation over $k$. What remains is the equality
\beq
\sum_{\Z\ni\alpha > -k_0}f_{\alpha\beta k_0}\:
\mu^{(m+\beta)\alpha}\:e_{m+\beta,k_0-1+\alpha}=0,
\;\;\; \Lall m,\beta\in\Z.
\eeq
Now the linear independence of all the summands implies that
\beq
f_{\alpha\beta k_0}=0,\;\;\; \forall\, \ha,\beta\in\Z,\, \alpha > -k_0\,.
\eeq
As \ha\ is always assumed to be bigger than $-k$, the condition $\alpha > -k_0$ is
automatically satisfied. This means that for all possible \ha\ and \hb\ we have
$f_{\alpha\beta k_0}=0$,
which contradicts the definition of $k_0$. Consequently, all coefficients
$f_{\alpha\beta k}$ must vanish, as desired. Thus $\rho(f)=0$ leads to the vanishing
of all coefficients $f_{\alpha\beta k}$ and $h_{\alpha\beta}$. Combining it with
$\rho'(f)=0$
yields
\beq
\sum_{\stackrel{\mbox{\scriptsize$\alpha,\beta,k\in\Z$}}{\beta\geq-n>-k}}
g_{\alpha\beta k}\,\mu^{-m\beta}\,
e_{m+\alpha,n+\beta}=0,\;\;\;\Lall m\in\Z, n\in\N.
\eeq
Much as before, we can deduce from here that also all coefficients $g_{\alpha\beta k}$
are zero. This proves the desired linear independence.

The second assertion of the lemma follows from the observation that in the above
reasoning the vanishing of all coefficients, which entails the equality $f=0$,
 is concluded solely from the formula
$(\rho\oplus\rho')(f)=0$.
\epf

This way we have shown the faithfulness of $\rho\oplus\rho'$ for any $0\leq p,q<1$. On the 
other hand, it is a general fact that, if $C$ is the universal $C^*$-algebra of a $*$-algebra
$\cO$, and the latter admits a faithful bounded $*$-representation, then $\cO\inc C$. 
Indeed, by the universality of $C$ we have the commutative diagram
\[
\begin{diagram}[height=10mm]
\cO&\rTo&B(\cH)\\
 &\rdTo&\uTo\\
 &&C
\end{diagram}.
\]
Hence, from the injectivity of the horizontal map (faithfulness of a representation), we can
infer the injectivity of the diagonal map, as needed. Consequently, we can claim:
\[
\cO(S^3_{pq\theta})\inc C(S^3_{pq\theta}),\;\;\;\forall\;0\leq p,q<1.
\] 

Any representation can be extended
 by continuity to the universal $C^*$-algebra, but in general there is no reason to expect that 
 the extension of a faithful representation remains faithful.  In particular, we cannot 
 automatically claim that the extension of $\rho\oplus\rho'$ is faithful. Here this difficulty
 is more pronounced as there is a problem with classifying irreducible representations,
 which is often a standard way for grasping the size of the universal $C^*$-algebra.
 (Recall that, in a sense, the irreducible representations play the role of points of the underlying
  quantum space.) Here a partial remedy to this problem is as follows.
\ble\label{faith}
Let $C_{\!A,p}$ and $C_{\!B,q}$ 
be the $C^*$-subalgebras of $C(S^3_{pq\theta})$, $0\leq p,q<1$,
generated by $a$ and~$b$, respectively. The extension of $\rho$ restricted to $C_{\!A,p}$
and the extension of $\rho'$ restricted to $C_{\!B,q}$ are faithful.
\ele
\bpf
Let us denote the aforementioned restrictions by $\rho_{\!A,p}$ and $\rho'_{\!B,q}\,$, respectively.
Due to the universality of the Topelitz algebra, the assignment $z\mapsto a$ defines
a $C^*$-homomorphism $\cT\cong C(D_p)\st{pr}{\ra} C_{\!A,p}$ (see (\ref{d})). Composing this map
with $\rho_{\!A,p}$ yields a representation of  $C(D_p)$ which splits into a direct sum of 
representations indexed by the integers corresponding to the first label of the basis
$\{e_{m,n}\}_{m\in\Z,\,n\in\N}$. It follows from Proposition~\ref{reps}(3) that the direct summand representation corresponding
to $m=0$ coincides with the extension %$\tilde\pi$ 
to $C(D_p)$ 
of the representation $\pi$ in (\ref{r}):
%Indeed, on the generator $z$ we have 
\[
\rho_{\!A,p}(pr(z))\,e_{0,n}=\sqrt{1-p^{n+1}}\;e_{0,n+1}\,.
\]
%which agrees with (\ref{r}). 
Since this extension is known to be faithful and $pr$ is evidently surjective, we can
conclude that $\rho_{\!A,p}$ is injective, as needed. The reasoning for $\rho'_{\!B,q}$ is identical.
\epf

\noindent
The argument proving the lemma shows also the injectivity of $pr$ and its $B,q$-counterpart.
This way, as both these maps are evidently surjective, we obtain as an immediate corollary that 
\[
C_{\!A,p}\cong\cT\cong B_{\!B,q}\,,~~~\forall\,p,q\in [0,1[\,.
\]

We are now ready to proceed to the key part of this section. It employs 
the just proved Lemma~\ref{faith} and is pivotal in the final $K$-theory computation.
\bth\label{00pq}
The $C^*$-algebra $C(S^3_{pq\theta})$, $0\leq p,q<1$, and $C(S^3_{00\theta})$
are isomorphic.
\ethe
\bpf
Assume first that $0<p,q<1$. Recall our convention that $a,b$ and $s,t$ stand for
 the generators of $C(S^3_{pq\theta})$ with $0<p,q<1$, and
$C(S^3_{00\theta})$, respectively. Since $s^*s=1=t^*t$ entails $\|s\|=1=\|t\|$, 
one can easily verify that the
infinite series
\[\label{fab}
\sum_{k=0}^\infty(\sqrt{1-p^{k+1}}-\sqrt{1-p^k})\,s^{k+1}{s^*}^k\;\mbox{ and }\;
\sum_{k=0}^\infty(\sqrt{1-q^{k+1}}-\sqrt{1-q^k})\,t^{k+1}{t^*}^k
\]
converge in norm to elements of $C(S^3_{00\theta})$. We denote them by $f(a)$ and
$f(b)$, respectively. On the other hand, (\ref{direl}) immediately implies the
positivity of the absolute values of generators: $|a|,|b|>0$. Hence $|a|,|b|$ are invertible,
and we can  define the following elements of $C(S^3_{pq\theta})$:
\[\label{gst}
g(s):=a|a|^{-1},\;\;\;
g(t):=b|b|^{-1}.
\]
The point is to extend the formulas (\ref{fab}) and (\ref{gst}) to mutually inverse 
$C^*$-homomorphisms. 
Our first step is to show the existence of $f$ and $g$. Due to the universality of both
$C^*$-algebras, it suffices to show that the elements in (\ref{fab}) and (\ref{gst})
satisfy appropriate relations. 

We begin with the sphere relation for the elements in (\ref{gst}). First, we have
\[
1-g(s)g(s)^*=1-a|a|^{-2}a^*
\]
 and
\[
(1-a|a|^{-2}a^*)(1-aa^*)=1-aa^*-a|a|^{-2}a^*+a(a^*a)^{-1}a^*aa^*=1-a|a|^{-2}a^*.
\]
In the same way, we obtain 
\[
1-g(t)g(t)^*=1-b|b|^{-2}b^*=(1-bb^*)(1-b|b|^{-2}b^*).
\]
Combining these facts, we compute
\[
(1-g(s)g(s)^*)(1-g(t)g(t)^*)
=
(1-a|a|^{-2}a^*)(1-aa^*)(1-bb^*)(1-b|b|^{-2}b^*)=0.
\]
Here the last step follows immediately from (\ref{sphere}). Next, we check the disc
relations:
\begin{align}
g(s)^*g(s)=|a|^{-1}a^*a|a|^{-1}=|a|^{-2}a^*a=1,~~~
g(t)^*g(t)=|b|^{-1}b^*b|b|^{-1}=|b|^{-2}b^*b=1.
\end{align}
Finally,  we employ (\ref{torel}) for the torus relations to obtain
\begin{align}
g(s)g(t)=a|a|^{-1}b|b|^{-1}
=e^{2\pi i\theta}ba|a|^{-1}|b|^{-1}
=e^{2\pi i\theta}b|b|^{-1}a|a|^{-1}
=e^{2\pi i\theta}g(t)g(s),
\end{align}
and, analogously, $g(s)g(t)^*=e^{-2\pi i\theta}g(t)^*g(s)$. Thus we have shown that
the formulas (\ref{gst}) define a $C^*$-homomorphism 
$g:C(S^3_{00\theta})\ra C(S^3_{pq\theta})$.

We now proceed to show that (\ref{fab}) defines~$f$. To make easier writing down computations, let us use the notation
\[
p_k:=\sqrt{1-p^{k+1}}-\sqrt{1-p^k},~~~q_k:=\sqrt{1-q^{k+1}}-\sqrt{1-q^k},~~~k\in\N.
\]
It is clear that $\sum_{k=0}^\infty p_k=1=\sum_{k=0}^\infty q_k$, and, consequently, that
$\sum_{k,l=0}^\infty p_kp_l=1=\sum_{k,l=0}^\infty q_kq_l$. Hence we can write
\begin{align}
(1-f(a)f(a)^*)(1-f(b)f(b)^*)
=
\!\!\!\sum_{k,l,m,n=0}^\infty \!\!\!p_kp_lq_mq_n
(1-s^{k+1}{s^*}^ks^l{s^*}^{l+1}) (1-t^{m+1}{t^*}^mt^n{t^*}^{n+1}).\nonumber
\end{align}
It follows from $s^*s=1=t^*t$ that the terms in parentheses are always of the form
$A_i$ and $B_j$, respectively. Thus the formula (\ref{akbl0}) entails the desired equality
\[
(1-f(a)f(a)^*)(1-f(b)f(b)^*)=0.
\]
Next, using $st=e^{2\pi i\theta}ts$ and $st^*=e^{-2\pi i\theta}t^*s$ we compute
\[
f(a)f(b)=\sum_{k,l=0}^\infty p_kq_ls^{k+1}{s^*}^kt^{l+1}{t^*}^l
%=\sum_{k,l=0}^\infty p_kq_lss^{k}{s^*}^ktt^{l}{t^*}^l
=e^{2\pi i\theta}\sum_{k,l=0}^\infty p_kq_lt^{l+1}{t^*}^ls^{k+1}{s^*}^k
=e^{2\pi i\theta}f(b)f(a)
\]
and, similarly,
\[
f(a)f(b)^*=e^{-2\pi i\theta}f(b)^*f(a).
\]
For the remaining  disc relations, we resort to
 the help of available representations.  Using definitions from Lemma~\ref{faith},
  it is straightforward to check
that 
\[\label{coi}
{\rho}_{\!A,0}(f(a))={\rho}_{\!A,p}(a),~~~ \rho'_{\!B,0}(f(b))=\rho'_{\!B,q}(b),~~~
{\rho}_{\!A,p}(g(s))={\rho}_{\!A,0}(s),~~~ \rho'_{\!B,q}(g(t))=\rho'_{\!B,0}(t).
\]
Therefore, since ${\rho}_{\!A,0}$ and $\rho'_{\!B,0}$ are injective and
 \[
 f(a)^*f(a)-pf(a)f(a)^*-1+p\in C_{\!A,0}\,,~~~
 f(b)^*f(b)-qf(b)f(b)^*-1+q\in C_{\!B,0}\,,
 \]
 the  disc relations (\ref{direl}) entail that
\[
 f(a)^*f(a)=pf(a)f(a)^*+1-p,~~~
 f(b)^*f(b)=qf(b)f(b)^*+1-q,
 \]
as needed. Hence the formulas (\ref{fab}) define a $C^*$-homomorphism 
$f:C(S^3_{pq\theta})\ra C(S^3_{00\theta})$.

It remains to show that $f$ and $g$ are mutually inverse. To this end, we shall again take
advantage of representations. First, note that the formulas (\ref{coi}) immediately imply
\[\label{coi2}
{\rho}_{\!A,p}={\rho}_{\!A,0}\ci f|_{C_{\!A,p}}\,,~~~ 
\rho'_{\!B,q}=\rho'_{\!B,0}\ci f|_{C_{\!B,q}}\,,~~~ 
{\rho}_{\!A,0}={\rho}_{\!A,p}\ci g|_{C_{\!A,0}}\,,~~~ 
\rho'_{\!B,0}=\rho'_{\!B,q}\ci g|_{C_{\!B,0}}\,.
\]
Hence we obtain
\begin{align}
{\rho}_{\!A,p}={\rho}_{\!A,p}\ci g|_{C_{\!A,0}}\ci f|_{C_{\!A,p}}\,,~~~
\rho'_{\!B,q}=\rho'_{\!B,q}\ci g|_{C_{\!B,0}}\ci f|_{C_{\!B,q}}\,,\\
{\rho}_{\!A,0}={\rho}_{\!A,0}\ci f|_{C_{\!A,p}}\ci g|_{C_{\!A,0}}\,,~~~
\rho'_{\!B,0}=\rho'_{\!B,0}\ci f|_{C_{\!B,q}}\ci g|_{C_{\!B,0}}\,.
\end{align}
 Consequently, 
the injectivity of all four representations used above (Lemma~\ref{faith})
implies that
\begin{align}
\id=g|_{C_{\!A,0}}\ci f|_{C_{\!A,p}}\,,~~~
\id= g|_{C_{\!B,0}}\ci f|_{C_{\!B,q}}\,,\\
\id=f|_{C_{\!A,p}}\ci g|_{C_{\!A,0}}\,,~~~
\id=f|_{C_{\!B,q}}\ci g|_{C_{\!B,0}}\,.
\end{align}
Therefore, as $a,b$ and $s,t$ generate $C(S^3_{pq\theta})$ and $C(S^3_{00\theta})$,
 respectively, we can conclude that $f$ and $g$ are mutually inverse, as desired. 
 This  can be summarized by the
commutative diagram
\[
\begin{diagram}[height=7mm,width=16mm]
& &C_{\!A,p}&\rTo^{\subseteq}&C(S^3_{pq\theta}) &\lTo^{\supseteq}&C_{\!B,q}&  \\
&\ldTo^{\rho_{\!A,p}}& & & & & & \rdTo^{\rho'_{\!B,q}}&\\
B(\cH)& &\dTo^{f|_{C_{\!A,p}}}\uTo_{g|_{C_{\!A,0}}} & 
& \dTo^{f} \uTo_{g}& 
& \uTo^{f|_{C_{\!B,q}}} \dTo_{g|_{C_{\!B,0}}}& &B(\cH)\\
&\luTo_{\rho_{\!A,0}}& & & & & & \ruTo_{\rho'_{\!B,0}}&\\
& &C_{\!A,0}&\rTo_{\subseteq}&C(S^3_{00\theta}) &\lTo_{\supseteq}&C_{\!B,0}&  \\
\end{diagram}.
\]

Finally, note that the formulas for $f$ and $g$ also make sense in the special cases when either $p$ or $q$ is equal to zero. Then the same argument shows that they give rise to mutually inverse $*$-homomorphisms. (For instance, if $p=0$, then the left column of the above  diagram is simply given by identity maps.)
\epf

Let us finish this section by discussing some special cases. 
Putting $\theta=0$ but keeping $0\leq p,q<1$ immediately recovers quantum-sphere
algebras in \cite{cm02,hms}.
On the other hand, since the norm
condition (\ref{norm}) is automatically satisfied for $0\leq p,q<1$, adding it to the definition
of  $C(S^3_{pq\theta})$, $0\leq p,q<1$, would leave it unchanged. Now,
however, taking for granted the theory of commutative $C^*$-algebras, it is clear from (\ref{class}) that setting in the thus defined  $C^*$-algebra $p=1=q$ and $\theta=0$
would yield $C(S^3)$. This is why we view our construction as a noncommutative deformation
of the topological 3-sphere. Furthermore, setting $p=1=q$ but keeping $\theta$ arbitrary
gives precisely the construction of the Matsumoto (constant $\theta$) $C^*$-algebra of a 
noncommutative 3-sphere \cite[p.334]{m-k91}.

\section{The fiber-product C*-algebra}
\setcounter{equation}{0}

The aim of this section is to prove that $C(S^3_{00\theta})$ is isomorphic to
a certain fiber product (pullback) of $C^*$-algebras, which will make it possible to compute
its $K$-groups using the Mayer-Vietoris exact sequence.

Let $\cT$ be the Toeplitz algebra, as defined in the preliminaries, 
and $z$ its generator. The crossed products ${\cal T}\rtimes_\theta\,\Z$
and ${\cal T}\rtimes_{-\theta}\,\Z$ with respect to the actions 
$\alpha^{\pm\theta}:\Z\times\cT\ra \cT$ given by
$
\ha^{\pm\theta}(n,z)=e^{\pm 2\pi in\theta}z
$
can be identified with the universal
$C^*$-algebras \cite[p.71]{b-b98} generated by $s_+,u$, respectively $t_-,v$, satisfying
\beq\label{o+}
s_+^*s_+=1=u^*u=uu^*,~~~s_+u=e^{2\pi i\theta}us_+\,,
\eeq
%respectively
\beq\label{o-}
t_-^*t_-=1=v^*v=vv^*,~~~t_-v=e^{-2\pi i\theta}vt_-\,.
\eeq
(Note that $s_+\mapsto t_-$, $u\mapsto v^*$ defines an isomorphism
$\cT\rtimes_\theta\Z\ra\cT\rtimes_{-\theta}\Z$.)

We interpret $\cT$ as the algebra of a quantum disc, and the just defined
crossed products as algebras of noncommutative solid tori. What we want to show is that
our noncommutative 3-sphere $S^3_{00\theta}$ is the quotient of the disjoint union of 
such two quantum
solid tori by an identification of their boundaries. 
The  common boundary of the two solid tori is the noncommutative 2-torus $T^2_\theta$
given by the (ir)rational rotation
algebra explained in preliminaries. 
The glueing
of the solid tori along their boundaries is with the exchange of fundamental cycles of
the boundaries, i.e., a cycle ``through a hole" is identified with a cycle ``around a hole".
This is reflected
in the exchange of the sign of $\theta$ in the construction of the quantum solid tori.

With this picture in mind, we have natural epimorphisms
\begin{align}
&{\pi}_1:{\cal T}\underset{\theta}{\rtimes}\Z\lra C(T^2_{\theta}),~~~
{\pi}_1(s_+):=x,~~~{\pi}_1(u):=y,
\\
&{\pi}_2:{\cal T}\underset{-\theta}{\rtimes}\Z\lra C(T^2_{\theta}), ~~~
{\pi}_2(t_-):=y,~~~{\pi}_2(v):=x.
\end{align}
Here $x$ and $y$ are the unitary generators of $C(T^2_{\theta})$ (see \ref{t2t}).
These two maps allow us to define the desired fiber product %algebra
of $\cT\rtimes_\theta\,\Z$ with $\cT\rtimes_{-\theta}\,\Z$ over $C(T^2_{\theta})$:
\[
\cT\underset{\theta}{\rtimes}\Z
\underset{C(T^2_{\theta})}{\oplus}
\cT\underset{-\theta}{\rtimes}\Z:=
\left\{(f_1,f_2)\in ({\cal T}\underset{\theta}{\rtimes}\Z)\oplus
({\cal T}\underset{-\theta}{\rtimes}\Z)~\Big|~{\pi}_1(f_1)={\pi}_2(f_2)\right\}.
\] 
Now we can define a homomorphism 
\[\label{h}
h:C(S^3_{00\theta})\lra 
\cT\underset{\theta}{\rtimes}\Z
\underset{C(T^2_{\theta})}{\oplus}
\cT\underset{-\theta}{\rtimes}\Z,~~~
h(s):=(s_+,v),~~~h(t):=(u,t_-).
\]
The situation is summarized in the following commutative diagram:
\[\label{sum}
\begin{diagram}[height=10mm,width=15mm]
\cT\underset{\theta}{\rtimes}\Z
\underset{C(T^2_{\theta})}{\oplus}
\cT\underset{-\theta}{\rtimes}\Z
&&\rTo^{\pr_2}&&
\cT\underset{-\theta}{\rtimes}\Z\\
&\luTo_{h~}&&\ruTo^{h_2~}&\\
\dTo^{\pr_1}&&
~~~~~~C(S^3_{00\theta})~~~~~~
&&\dTo_{\pi_2}\\
&\ldTo_{~h_1}&&\rdTo^{~\phi_c}&\\
\cT\underset{\theta}{\rtimes}\Z&&\rTo_{\pi_1}&&C(T^2_\theta)
\end{diagram}.
\]
Here $\pr_1$ and $\pr_2$ are the natural projections, and 
$h_1$, $h_2$, $\phi_c$ are the appropriate composite maps.
Note that all these homomorphisms introduced above are well-defined 
on the $C^*$-level because
their domains of definition are universal $C^*$-algebras.

Our goal is to show that $h$ is an isomorphism. To this end, we need to introduce
some more $C^*$-homomorphisms whose existence is guaranteed by the following
general result:
\ble\label{gen}
Let $C$ be a (not necessarily unital) 
$C^*$-algebra, $w$ the unitary generator of  the $C^*$-algebra
 of continuous functions on $S^1$, and $E_{ij}$ the matrices with zero everywhere except for
 the $ij$-place where there is 1 (matrix units). The assignment $E_{ij}\ot 1\mapsto e_{ij}$,
$E_{ij}\ot w\mapsto w_{ij}$, $e_{ij},w_{ij}\in C$, defines a homomorphism
$\cK\ot C(S^1)\ra C$ if and only if the following conditions are satisfied
for any $N\in\N$:
\vspace*{-1mm}\begin{itemize}
\item[(1)] the $e_{ij}$ fulfill the relations of $E_{ij}$, i.e., 
$e_{ij}e_{kl}=\delta_{jk}e_{il}$, $e_{ij}^*=e_{ji}$, 
$\forall~ 0\leq i,j,k,l\leq N$;
\item[(2)] $W:=\sum_{i=0}^Nw_{ii}$ is a partial unitary, i.e., $W^*W=WW^*$
is a projector;
\item[(3)] $[W,e_{ij}]=0$, $e_{ij}W=w_{ij}$, $\forall~ 0\leq i,j\leq N$.
\end{itemize}
\ele
\bpf
Note first that for a given homomorphism
$\cK\ot C(S^1)\ra C$, the values on $E_{ij}\ot 1$ and $E_{ij}\ot w$ obviously
have to satisfy the above conditions. For the proof the other way around, choose
an arbitrary $N\in \N$ and consider  $M_{N+1}(\C)\ot C(S^1)$. To show that there exists
a homomorphism from this $C^*$-algebra to $C$, by the universality of the tensor 
product (\cite[Proposition~4.7]{t-m79}), it suffices to show that there exist
 homomorphisms $M_{N+1}(\C)\ra C$ and $C(S^1)\ra C$ whose ranges commute. Since 
 $E_{ij}$ form a basis of $M_{N+1}(\C)$, the condition (1) guarantees that the assignment
 $E_{ij}\mapsto e_{ij}$ defines a $C^*$-algebra homomorphism. On the other hand,
 for any bounded normal operator $U$ such that $(U^*U)^2=U^*U$, we have  $U=U(U^*U)$.
 Indeed, as $|U|:=\sqrt{U^*U}$ can be approximated by polynomials divisible by $U^*U$,
 the desired formula follows from the polar decomposition: 
 \[
 U(U^*U)=\mbox{Phase}(U)\,|U|\,U^*U
 =\mbox{Phase}(U)\lim_{n\ra\infty}(p_n(U^*U)U^*U)\,U^*U
 =U.
\]
Therefore, since every $C^*$-algebra admits a faithful representation, we can conclude from (2) that $W=W(W^*W)$. Now it follows from the normality of $W$ that $W^*W$ is the identity
of the $C^*$-subalgebra generated by $W$. Hence the assignment $w\mapsto W$ defines
a homomorphism $C(S^1)\ra C$. By (3), the ranges of the two just constructed 
homomorphisms commute. Thus, for any $N\in\N$, we have a homomorphism
$M_{N+1}(\C)\ot C(S^1)\ra C$. These homomorphisms define a $*$-algebra homomorphism
$\bigcup_{N\in\N}M_{N+1}(\C)\ot C(S^1)\ra C$. Since the value of this homomorphism on
any element is given by applying a $C^*$-homomorphism, its norm is preserved or 
decreased. Consequently, the $*$-algebra homomorphism is norm-continuous and it
extends to a homomorphism:
\[
\cK\ot C(S^1)=\overline{\bigcup_{N\in\N}M_{N+1}(\C)\ot C(S^1)}\overset{F}{\lra} C.
\]
To end the proof, it suffices to note that $F(E_{ij}\ot 1)=e_{ij}$ and that the condition
(3) entails
$F(E_{ij}\ot w)=w_{ij}$.
\epf

We are now ready to prove the main statement of this section:
\bth\label{fibreuni}
The homomorphism (\ref{h}) is an isomorphism of the universal $C^*$-algebra $C(S^3_{00\theta})$ with the fiber-product $C^*$-algebra 
$\cT\rtimes_{\theta}\,\Z
\;\oplus_{C(T^2_{\theta})}\;
\cT\rtimes_{-\theta}\,\Z$.
\ethe
\bpf 
Our strategy is to construct a commutative diagram of two short exact sequences and then conclude the assertion of the theorem by  the Five Isomorhisms Lemma.
To begin with, one can verify  with the help of Lemma~\ref{gen} that the formulas
\begin{align}
&(E_{ij}\ot w,0)\longmapsto s^i(1-ss^*)t{s^*}^j,~~~
(0,w\ot E_{ij})\longmapsto t^i(1-tt^*)s{t^*}^j,~~~ \mbox{and}
\\
&(E_{ij}\ot w,0)\longmapsto (s_+^i(1-s_+s_+^*)u{s_+^*}^j,0),~~~
(0,w\ot E_{ij})\longmapsto (0,t_-^i(1-t_-t_-^*)v{t_-^*}^j),
\end{align}
where $E_{ij}$ are the matrix units and $w$ is the unitary generator of
$C(S^1)$,
define homomorphisms
\begin{align}
&j_c:(\cK\ot C(S^1))\oplus(C(S^1)\ot \cK)\lra C(S^3_{00\theta})~~~ \mbox{and}
\\
&j_d:(\cK\ot C(S^1))\oplus(C(S^1)\ot\cK)\lra 
\cT\rtimes_\theta\Z\!\underset{C(T^2_\theta)}{\oplus}\!\cT\rtimes_{-\theta}\Z,
\end{align}
 respectively. Observe also that the above formulas entail
\begin{align}
&j_c(E_{ij}\ot 1,0)= s^i(1-ss^*){s^*}^j,~~~
j_c(0,1\ot E_{ij})= t^i(1-tt^*){t^*}^j,
\label{jc}\\
&j_d(E_{ij}\ot 1,0)= (s_+^i(1-s_+s_+^*){s_+^*}^j,0),~~~
j_d(0,1\ot E_{ij})= (0,t_-^i(1-t_-t_-^*){t_-^*}^j).
\end{align}
Next, define $\phi_d$ as the composition $\pi_1\ci\pr_1=\pi_2\ci\pr_2$ and recall
that $\phi_c:=\phi_d\ci h$ (see (\ref{sum})). One can immediately check the commutativity
of the diagram:
\[\label{5}
\begin{diagram}[height=10mm]
0 & \rTo & (\cK\ot C(S^1))\oplus(C(S^1)\ot\cK) & \rTo{~~j_c} &
 C(S^3_{00\theta})
& \rTo{\phi_c} & C(T^2_\theta)& \rTo & 0\,\phantom{.}\\
 && \dTo_{\id} && \dTo_{h} && \dTo_{\id} && \\
0 & \rTo &(\cK\ot C(S^1))\oplus(C(S^1)\ot \cK) & \rTo{j_d} & \cT\underset{\theta}{\rtimes}\Z
\underset{C(T^2_{\theta})}{\oplus}
\cT\underset{-\theta}{\rtimes}\Z & 
\rTo{\phi_d}
 & C(T^2_\theta)& \rTo & 0\,. 
\end{diagram}
\]

A key step is to prove the exactness of the rows.
Consider first the upper row.
In order to show the injectivity of $j_c$,  we will show that
$(\rho\oplus\rho')\ci j_c$ is injective. (Here $\rho$ and $\rho'$ are representations
from Proposition~\ref{reps}.) Since $(\rho\ci j_c)(0,T')=0=(\rho'\ci j_c)(T,0)$,
 this is the case if and only if the homomorphisms
 \begin{align}
&\psi:\cK\ot C(S^1)\lra B(\cH), ~~~ \psi(T):=(\rho\ci j_c)(T,0), 
\\
&\psi':C(S^1)\ot\cK\lra B(\cH),~~~\psi'(T'):=(\rho'\ci j_c)(0,T'),
\end{align}
 are both injective. 
Let us check the injectivity of $\psi$.
%Assume that
%$\ker((\rho\oplus\rho')\ci j_c=\ker(\rho\ci j_c)\cap\ker(\rho'\ci j_c)\neq 0$.
%Then it follows from $(0,\cK\ot C(S^1))\subset \ker\rho\ci j_c$ and
%$(\cK\ot C(S^1),0)\subset \ker\rho'\ci j_c$
%that there must exist ideals $J, J'\subset \cK\ot C(S^1)$ with at least
%one of them being different from zero, such that
%$(J,\cK\ot C(S^1))=\ker\rho\ci j_c$ and $(\cK\ot C(S^1),J')=\ker\rho'\ci j_c$.
%Assume that $J\neq 0$. 
Clearly, $\ker\psi$ is a closed ideal in $\cK\ot C(S^1)$. Hence it is of the form
$\cK\ot I$, with $I$ an ideal in $C(S^1)$. 
In particular, for any $f\in I$, we have $E_{00}\ot f\in\ker\psi$. Now one checks
that 
\[
\psi(E_{00}\ot w)\,e_{m,n}=\delta_{0n}\,e_{m+1,0}\,,~~~\mbox{i.e.,}
~~~\psi(E_{00}\ot w)=\sigma_d\ot E_{00}\,,
\]
 where $\sigma_d$ is the
standard two-sided shift. Therefore, since the assignment $w\mapsto \sigma_d$ 
defines a faithful
representation $\rho_d$ of $C(S^1)$, for any $f\in I$ we have:
\[
0=\psi(E_{00}\ot f)=\rho_d(f)\ot E_{00}\;\;\;\imp\;\;\; f=0.
\] 
Consequently, $\ker\psi=0$. One argues analogously for $\psi'$.
%if $J'\neq 0$ is assumed.
Thus, the injectivity of $j_c$ is shown.
For the middle exactness, note that $\phi_c\ci j_c=0$, so that $\phi_c$ induces
 the map
 \[ \label{map}
C(S^3_{00\theta})\;\big/\;j_c\llp(\cK\ot C(S^1))\oplus(\cK\ot C(S^1))\lrp\lra C(T^2_\theta).
\]
On the other hand, it follows from (\ref{o+}--\ref{o-}) and (\ref{jc}) that
the classes of $s$ and $t$ in the quotient algebra satisfy the relations (\ref{t2t})
of the noncommutative torus. Therefore, due to the universality of the $C^*$-algebra
$C(T^2_\theta)$, there is a natural homomorphism in the
other direction, which is evidently the inverse of (\ref{map}).
This finishes the proof of exactness of the upper sequence.

The proof for the lower sequence proceeds along the same lines. This time, instead of $\rho$ and $\rho'$ one uses 
$\bar{\rho}\ci\pr_1$ and $\bar{\rho}'\ci\pr_2$, where $\bar{\rho}$ and $\bar{\rho}'$ 
 are representations  of 
$\cT\rtimes_\theta\Z$ and  $\cT\rtimes_{-\theta}\Z$, respectively,
defined by
\begin{align}
&\bar{\rho}(s_+)\,e_{m,n}=\mu^m \,e_{m,n+1},~~~\bar{\rho}(u)\,e_{m,n}=e_{m+1,n},
\\
&\bar{\rho}'(t_-)\,e_{m,n}=\mu^{-m}\, e_{m,n+1},~~~\bar{\rho}'(v)\,e_{m,n}=e_{m+1,n}.
\end{align}
These representations are related to $\rho$ and $\rho'$ by the maps $h_1$ and $h_2$
of the diagram (\ref{sum}):
 $\rho=\bar{\rho}\ci h_1$ and $\rho'=\bar{\rho}'\ci h_2$.
To finish the exactness proof, we identify
\[
\left(\cT\underset{\theta}{\rtimes}\Z
\underset{C(T^2_{\theta})}{\oplus}
\cT\underset{-\theta}{\rtimes}\Z\right)\;\Big/\;j_d\llp(\cK\ot C(S^1))\oplus(\cK\ot C(S^1))\lrp
\]
with  $C(T^2_{\theta})$ by analysing the classes of generators $(s_+,v)$ and $(u,t_-)$.
Now, the proof of the theorem is completed by applying the Five Isomorphisms Lemma
to the diagram (\ref{5}).
\epf

One can infer from the above proof that there is a short exact sequence
\[
0\lra\cK\ot C(S^1)\lra\cT\underset{\theta}{\rtimes}\Z\lra 
C(S^1)\underset{\theta}{\rtimes}\Z =C(T^2_{\theta})\lra 0.
\]
We can interpret it as the decomposition of our quantum solid torus into its inside and
the boundary. An interesting phenomenon manifests itself in this decomposition. It makes
this quantum solid torus significantly different from the one considered by Matsumoto 
\cite{m-k91},
notably $C(D^2)\rtimes_{\theta}\Z$. The inside of the latter torus is composed out of the 
nested noncommutative tori shrinking to a circle, so that the action of $\Z$ is not restricted
to the boundary of the disc $D^2$. On the contrary, in our case there is no parameter
$\theta$ in the ideal representing the inside of the quantum solid torus.
 In this sense, the noncommutative
deformation of the algebra of the interior of $D^2$ into $\cK$ pushed out the 
$\theta$-deformation. 

\section{The K-groups}
\setcounter{equation}{0}

In what follows, we compute the $K$-theory of the fiber-product $C^*$-algebra 
${\mathcal T}\!{\rtimes}_ {\theta}\,{\mathbb Z}
\,{\oplus}_{C(T^2_\theta)}\,
{\mathcal T}\!{\rtimes}_{-\theta}\, {\mathbb Z}$ from
the previous section. This, combined with identifications of this fiber product %$C^*$-algebra
with the universal $C^*$-algebras of Section~2, yields the main result of our paper, which is the computation of the $K$-groups of the Heegaard-type quantum 3-spheres.

As the first step, let us  determine
the $K$-theory of the crossed product 
$\cT{\rtimes}_ {\theta}\,\Z\cong\cT{\rtimes}_ {-\theta}\,\Z$. We know that
$K_0(\cT)\cong\Z$ (generated by the class of $1\in\cT$) and $K_1(\cT)=0$ 
\cite[p.179]{w-ne93}.
Therefore, since the action of $\Z$ on $1\in\cT$ is trivial, 
the Pimsner-Voiculescu exact sequence \cite[Theorem~10.2.1]{b-b98} reduces to
\[
0\lra K_1 \left(\cT\!\underset{\theta}{\rtimes}\Z\right)\lra\Z\overset{0}{\lra}\Z
\lra K_0 \left(\cT\!\underset{\theta}{\rtimes}\Z\right) \lra 0.
\]
This immediately leads to
\[\label{kpro}
K_j \left(\cT\!\underset{\pm\theta}{\rtimes}\Z\right)\cong\Z,~~~ j\in\{0,1\},
\]
and gives explicitly 2 out of the 6 terms of
 the Mayer-Vietoris exact sequence \cite[Example~21.1.2(a) and Theorem~21.2.2]{b-b98}:
 
\beq\label{qmv}
\begin{diagram}[height=12mm,width=12mm]
K_0\left(\cT\!\underset{\theta}{\rtimes}\Z
\underset{C(T^2_\theta)}{\oplus} 
\cT\!\underset{-\theta}{\rtimes}\Z\right) 
&  \rTo 
& K_0\left( \cT\!\underset{\theta}{\rtimes}\Z\right) \oplus
 K_0\left(\cT\!\underset{-\theta}{\rtimes}\Z\right)
&\rTo  
&  K_0(C(T^2_\theta)) \\
 \uTo &  &  & &\dTo \\
K_1(C(T^2_\theta))& \lTo     
&K_1\left( \cT\!\underset{\theta}{\rtimes}\Z\right) \oplus
 K_1\left(\cT\!\underset{-\theta}{\rtimes}\Z\right) &\lTo 
&
 K_1\left(\cT\!\underset{\theta}{\rtimes}\Z
\underset{C(T^2_\theta)}{\oplus} 
\cT\!\underset{-\theta}{\rtimes}\Z\right) 
\end{diagram}.%~~~~~~~~~~
\eeq
For further unravelling of this sequence recall that
$K_0 (C(T^2_\theta))\cong {\mathbb Z} \oplus {\mathbb Z}$ \cite[Section~12.3]{w-ne93}. 
Here the first ${\mathbb Z}$ is generated by the unit of the $C^*$-algebra 
$C(T^2_\theta)$. The second ${\mathbb Z}$ is generated by the ``Hopf line bundle'' on the noncommutative torus $T^2_\theta$. Since 
$K_0 ({\mathcal T}{\rtimes}_ {\pm\theta}\,{\mathbb Z}) \cong{\mathbb Z}$ and is generated by the unit of  ${\mathcal T}{\rtimes}_ {\theta}\,{\mathbb Z}$,  
the upper right horizontal arrow becomes
\[
{\mathbb Z} \oplus {\mathbb Z}\ni (m,n) \longmapsto (m-n , 0)\in
{\mathbb Z} \oplus {\mathbb Z}.
\]
On the other hand,
$K_1 (C(T^2_\theta))\cong  {\mathbb Z} \oplus {\mathbb Z}$,
 where the generators are the classes of the two unitary generators of $C(T^2_\theta)$.
Therefore, since 
$K_1 ({\mathcal T}{\rtimes}_ {\pm\theta}\,{\mathbb Z}) \cong {\mathbb Z}$
 is generated by the unitary coming from the action of ${\mathbb Z}$,
  the lower left horizontal arrow becomes
\[
{\mathbb Z} \oplus {\mathbb Z}\ni (m,n) \longmapsto (m,-n)\in
{\mathbb Z} \oplus {\mathbb Z}.
\]
It is evidently an isomorphism, whence the preceding and following arrows are zero maps.
Summarizing, with the shorthand notation $G_0$ and $G_1$ for the left-top and right-bottom corners of (\ref{qmv}), respectively, we have obtained:
\beq\label{cqmv}
\begin{diagram}[height=10mm]
G_0
&  \rTo 
& \Z \oplus\Z
&\rTo{(m,n)\mapsto(m-n,0)}  
&  \Z \oplus\Z \\
 \uTo{0} &  &  & &\dTo \\
\Z \oplus\Z& \lTo{(m,-n)\mapsfrom(m,n)}    
&\Z \oplus\Z&\lTo{0} 
&
G_1
\end{diagram}.
\eeq
The exactness of this sequence entails that
\begin{align}
G_0&\cong\ker\{\Z \oplus\Z\ni(m,n)\longmapsto(m-n,0)\in \Z \oplus\Z\}\cong\Z,
\label{g0}\\
G_1&\cong\coker\{\Z \oplus\Z\ni(m,n)\longmapsto(m-n,0)\in \Z \oplus\Z\}\cong\Z
\label{g1}.
\end{align}

Finally, combining (\ref{g0}--\ref{g1}) with Theorem~\ref{fibreuni} and Theorem~\ref{00pq} 
yields our main claim establishing the $K$-groups of the universal $C^*$-algebra of
Definition~\ref{cdef}:
\bth\label{main}
$K_j (C(S^3_{pq\theta}))\cong\Z$, $j\in\{0,1\}$, $0\leq p,q,\theta<1$.
\ethe

%\small\footnotesize
{\bf Acknowledgments.}
Much of the joint work on this paper took place in
the Institut Mittag-Leffler (IML)
 in Djursholm and the Institut des Hautes \'{E}tudes Scientifiques (IHES)
in Bures-sur-Yvette.
It is a pleasure to thank these institutes for their wonderful 
hospitality and support. Here special thanks are due to
C.~ Cheikhchoukh, F.~Schmit and M.C.~Vergne of IHES for their help with
preparing the manuscript.
P.M.H.\ is also grateful to S.L.~Woronowicz for some $C^*$-algebraic inspiration,
and to M.~Rieffel for discussions.
This work was partially supported by
the  KBN grant 2 P03A 013 24 (P.M.H., R.M. and W.S.), IML and IHES 
visiting stipends (P.B. and P.M.H.),
and  the Graduiertenkolleg ``Quantenfeldtheorie" at 
Universit\"at Leipzig (P.M.H.). 
%The authors are grateful  to 

\end{document}